\documentclass[11pt]{amsart}
\ifdefined\XeTeXversion
  \PassOptionsToPackage{xetex}{hyperref}
\else
  \ifdefined\pdfoutput\pdfoutput=1\fi
\fi

\usepackage[usenames,dvipsnames]{xcolor}   %

\usepackage[utf8]{inputenc}
\usepackage[T1]{fontenc}

\usepackage[letterpaper,margin=1in]{geometry}

\usepackage{lmodern}
\usepackage{microtype}
\usepackage{etoolbox}                 %
\usepackage{lettrine}

\usepackage{mathtools}
\usepackage{amsmath,amssymb,amsfonts,amsxtra,bm}
\usepackage{amsthm}
\usepackage{mathrsfs}
\usepackage{latexsym}
\usepackage{nicefrac}       %

\usepackage{booktabs}
\usepackage{makecell,multirow}       %
\usepackage{colortbl}
\usepackage{graphicx}
\usepackage{epstopdf}
\usepackage{float}
\floatstyle{plaintop}
\restylefloat{table}
\usepackage{footnote, tablefootnote}

\usepackage{enumitem}
\usepackage{algorithm}
\usepackage{algorithmic}
\usepackage{comment}
\usepackage{xspace}
\usepackage{url}
\usepackage[most]{tcolorbox}

\usepackage[round]{natbib}

\usepackage{hyperref}

\numberwithin{equation}{section}

\makeatletter

\definecolor{darkblue}{rgb}{0.0,0.0,0.65}
\definecolor{darkred}{rgb}{0.68,0.05,0.0}
\definecolor{darkgreen}{rgb}{0.0,0.29,0.29}
\definecolor{darkpurple}{rgb}{0.47,0.09,0.29}

\patchcmd{\section}{\scshape}{\fontfamily{ppl}\bfseries\scshape\selectfont}{}{%
  \PackageWarning{aiproofs}{Could not bold the section font; amsart may have changed}}

\newcommand{\aiauthor}{%
  \author[Suvrit Sra]{Suvrit Sra\textsuperscript{\normalfont*}}%
  \thanks{*\,TU Munich, Dept.~of Mathematics, School of CIT; Garching, Germany.
    Email: \texttt{s.sra@tum.de}.}}

\newcommand{\arxivlink}[1]{%
  \thanks{Also on arXiv: \url{https://arxiv.org/abs/#1}.}}
\newcommand{\preprintsfooter}[2]{%
  \def\ps@firstpage{\ps@plain
    \def\@oddfoot{\hfil\parbox[t]{\textwidth}{\centering\normalfont\scriptsize
      First public release: #2, as v1 in the author's preprint repository,
      \url{https://github.com/suvrit/preprints/tree/main/preprints/#1};
      every version is kept there with its timestamp.\\[3pt]\thepage}\hfil
      \global\topskip\normaltopskip}%
    \let\@evenfoot\@oddfoot
    \def\@oddhead{\@serieslogo\hss}%
    \let\@evenhead\@oddhead}}

\makeatother

\hypersetup{colorlinks=true,linkcolor=red!60!black,citecolor=green!48!black,urlcolor=blue!60!black}
\setlist{nosep,leftmargin=*}

\newcommand{\contriblabel}[1]{{\color{darkblue}\bfseries(#1)}}
\newlist{contributions}{enumerate}{1}
\setlist[contributions]{label=\contriblabel{\arabic*},ref=(\arabic*),
  leftmargin=2.4em,labelsep=0.5em,itemsep=1pt,topsep=2pt}

\newenvironment{aiusage}{\section*{AI usage}}{}

\newcommand{\aiacknowledgements}[1][]{%
  \section*{Acknowledgements}
  \if\relax\detokenize{#1}\relax\else #1\ \fi
  The author (SS) thanks the Alexander von Humboldt (AvH) Foundation for
  their generous support via a Humboldt Professorship in AI.%
}

\newcommand{\C}{\mathbb{C}}

\newcommand{\etal}{et~al.\@\xspace}

\newcommand{\nlsum}{\sum\nolimits}
\newcommand{\nlprod}{\prod\nolimits}
\newcommand{\ip}[2]{\langle {#1},\, {#2} \rangle}

\newcommand{\norm}[1]{\|{#1}\|}

\DeclareMathOperator{\diag}{diag}

\DeclareMathOperator{\per}{per}

\DeclareMathOperator{\Imm}{Imm}

\newtheorem{theorem}{Theorem}[section]

\newtheorem{prop}[theorem]{Proposition}
\newtheorem{lemma}[theorem]{Lemma}
\newtheorem{corr}[theorem]{Corollary}
\newtheorem{rmk}[theorem]{Remark}
\newtheorem{example}[theorem]{Example}

\theoremstyle{definition}

\definecolor{cdarkred}{RGB}{180,70,30} %
\definecolor{cdarkgreen}{RGB}{0,130,0} %
\definecolor{azure}{RGB}{240,255,255}
\definecolor{cmidgreen}{RGB}{201,219,216}
\definecolor{clightgreen}{RGB}{200,255,200}
\definecolor{clightgray}{RGB}{240,240,240}
\definecolor{cmygreen}{RGB}{10,120,0}
\definecolor{cmyblue}{RGB}{50,50,255}
\definecolor{sblue}{RGB}{100,100,255}
\definecolor{cmypurple}{RGB}{150,50,180}
\definecolor{cmyred}{RGB}{238, 64,  0}
\definecolor{cdarkblue}{RGB}{30,75,170}
\definecolor{csharpred}{RGB}{255,0,0}

\tcbset{dashedrounded/.style={no borderline,borderline={0.5pt}{0pt}{dashed},arc=4pt}}
\tcbset{supporting/.style={colback=gray!8!white}}
\newtcolorbox{dashedbox}[1][]{
  before={\par\medskip\noindent},
  after={\par\medskip},
  colback=yellow!6!white,
  boxrule=0pt,
  enhanced jigsaw,
  top=2pt,
  left=2pt,
  right=2pt,
  bottom=0pt,
  dashedrounded,
  #1
}

\newtcbox{\tcbmathbox}[1][]{
  nobeforeafter,
  colback=yellow!30!white,
  boxrule=0pt,
  enhanced jigsaw,
  top=0pt,
  left=1pt,
  right=1pt,
  bottom=0pt,
  borderline horizontal={0.5pt}{0pt}{dashed},
  borderline vertical={0.5pt}{0pt}{dashed},
  #1
}

\hypersetup{pdftitle={The generalized Lax conjecture holds for mixed determinants},pdfauthor={Suvrit Sra}}
\allowdisplaybreaks[2]

\newcommand{\R}{\mathbb R}
\newcommand{\Q}{\mathbb Q}
\newcommand{\F}{\mathbb F}
\newcommand{\HH}{\mathcal H}
\newcommand{\eps}{\varepsilon}
\DeclareMathOperator{\tr}{tr}

\DeclareMathOperator{\Span}{span}
\DeclareMathOperator{\op}{op}
\newcommand{\cone}{\Lambda_+}
\newcommand{\ocone}{\Lambda_{++}}
\newcommand{\one}{\mathbf 1}
\newcommand{\pdm}{p_{d,m}}
\newcommand{\Ldm}{L_{d,m}}


\begin{document}
\title[Lax conjecture for mixed determinants]{\makebox[\textwidth]{\textls[-20]{The generalized Lax conjecture holds for mixed determinants}}}
\aiauthor
\begin{abstract}
  We prove that the mixed determinant $\eta(\cdot)$ of \citet{branden08} has a spectrahedral
  hyperbolicity cone. More precisely, for Hermitian matrices $A_1,\ldots,A_m$ of order $d$, the
  hyperbolicity cone of $\eta(A_1,\ldots,A_m)$ is cut out by $\nlsum_j\Phi_j(A_j)\succeq0$, where
  $\Phi_1,\ldots,\Phi_m$ are explicit unital completely positive maps. We argue by induction on $d$,
  adding one row and one column to every argument and expressing the logarithmic derivatives of $\eta$
  as extremal eigenvalues. By compressing along polynomial boundary kernels, we reduce the order of the
  pencil $\nlsum_j\Phi_j(A_j)$ to $\exp(O(md^3\log m))$; for one fixed positive semidefinite $n\times n$
  matrix $A$ a separate pencil yields the order $\exp(O(n^2\log n))$. We also answer four problems of
  Borcea and Br\"and\'en
  negatively. %
\end{abstract}
\thanks{First public release: 2 November 2026, as v1 in the author's preprint repository,
  \url{https://github.com/suvrit/preprints/tree/main/preprints/eta}, where every version is kept with
  its timestamp; also on Hexagon, \url{https://hexagonmath.org/abs/2611.00117}, and on the author's
  page, \url{https://suvrit.github.io}.}
\maketitle
\addtocontents{toc}{\protect\setcounter{tocdepth}{-1}}

\vspace*{-6pt}
\section{Introduction}\label{sec:intro}
A homogeneous polynomial $p$ on a real vector space is \emph{hyperbolic} with respect to direction $e$
if $p(e)\ne0$ and for every $x$ the polynomial $t\mapsto p(x+te)$ has only real zeros. Write
$\ocone(p,e)$ for the connected component of $\{p\ne0\}$ containing $e$, and $\cone(p,e)$ for its
closure, the hyperbolicity cone. A cone is \emph{spectrahedral} if it has the form $\{x:L(x)\succeq0\}$
for a homogeneous real symmetric linear pencil $L$, so that optimization over it is a semidefinite
program. If $L(e)\succ0$, then $\det L$ is hyperbolic with respect to $e$ and its hyperbolicity cone
is $\{x:L(x)\succeq0\}$. Lax conjectured that every hyperbolic polynomial in three variables arises in
this way, up to a constant factor~\citep{lax1958differential}; the conjecture is
true~\citep{helton2007linear,lewis2005lax}, but in more variables the statement fails, even for powers
of the polynomial~\citep[Theorem~3.3]{branden2011obstructions}.

However, for optimization only the cone matters, and a pencil can cut out the cone even when its
determinant carries extra factors. Thus, the generalized Lax conjecture asks only that every
hyperbolicity cone be spectrahedral.  \citet{branden2014hyperbolicity} proved it for the elementary
symmetric polynomials, and later papers confirmed it in further
cases~\citep{kummer2016vamos,saunderson2018spectrahedral,amini2019spectrahedrality,kummer2021spectral,kummer2026strictly}.
Unfortunately, the general conjecture fails, as shown by OpenAI's recent counterexample~\citep{openai2026nonspectrahedral}. The earlier successes still raise a positive question:
\emph{which families of hyperbolic polynomials do have spectrahedral cones?}

\vskip5pt
\noindent\textbf{Our choice.} We study mixed determinants, a remarkable family of hyperbolic polynomials introduced by~\citet{branden08}. Let $\F\in\{\R,\C\}$, and let $\HH_d(\F)$ denote the real vector space of $d\times d$ real symmetric matrices when $\F=\R$, or of $d\times d$ complex Hermitian matrices when $\F=\C$. For $A\in\HH_d(\F)$ and $S\subseteq[d]=\{1,\ldots,d\}$, write $A[S]$ for the principal submatrix of $A$ indexed by $S$, with $\det A[\varnothing]=1$. The mixed determinant of $m$ arguments of matrix order $d$ is
\begin{equation}
\label{eq:1}
 \pdm(A_1,\ldots,A_m)
 =\eta(A_1,\ldots,A_m)
 :=\sum_{S_1\sqcup\cdots\sqcup S_m=[d]}
       \nlprod_{j=1}^m\det A_j[S_j].
\end{equation}
Note that the mixed determinant is \emph{not} the mixed discriminant, a different object. We retain the notation $\eta$ of Borcea and Br\"and\'en for specialized arguments, but write $\pdm$ when the sizes $d$ and $m$ matter. For $A\succeq0$ and $B$ Hermitian, $\eta(zA,-B)$ is the polynomial arising in Johnson's
conjectures~\citep{johnson1987permanent,johnson1989characteristic}, proved by~\citet{branden08}.

The mixed determinant~\eqref{eq:1} is a real homogeneous polynomial of degree $d$ in the independent real coordinates of its arguments. Set $e_{d,m}=(I_d,0,\ldots,0)$. Then $\pdm(e_{d,m})=1$, and by Theorems~2.6 and~2.8 of \citet{branden08} (Lemma~\ref{lem:interlacing} below), $\pdm$ is hyperbolic with respect to $e_{d,m}$, with hyperbolicity cone $\mathcal K_{d,m}=\cone(\pdm,e_{d,m})$. For a rank-one matrix $A=vv^*$, every principal minor of order at least two vanishes, and
\begin{equation}\label{eq:rankone-intro}
 \eta(z_1A,\ldots,z_mA)
 =d!\left(\nlprod_{i=1}^d|v_i|^2\right)e_d(z_1,\ldots,z_m),
\end{equation}
where $e_d$ is the elementary symmetric polynomial of degree $d$. Given that cones of elementary symmetric polynomials are spectrahedral, observation~\eqref{eq:rankone-intro} raises the natural question whether for higher-rank $A$ the associated cone of $\eta$ is also spectrahedral.

The question reaches beyond $\eta$: Borcea and Br\"and\'en asked whether every homogeneous
\emph{symmetric} real stable polynomial is of the form $\eta(z_1A,\ldots,z_mA)$ with
$A\succeq0$~\citep[Problem~35]{aim2007polyaschurlax}, so a positive answer would have made the
hyperbolicity cones of all such polynomials spectrahedral (Appendix~\ref{sec:aim} shows that the answer
is negative). Our main result (Theorem~\ref{thm:main}) represents the cone $\mathcal K_{d,m}$ by a
pencil whose coefficient maps depend only on $d$ and $m$, so a single pencil serves every choice of the
matrices; we call $\pdm$, its cone $\mathcal K_{d,m}$, and this pencil \emph{universal}.

\begin{dashedbox}
\begin{theorem}\label{thm:main}
For every $d,m\geq1$, there are integers $N=N_{d,m}$ and $r=r_{d,m}$ and fixed matrices $K_{j,k}\in\{0,1\}^{N\times d}$, for $1\leq j\leq m$ and $1\leq k\leq r$, such that the maps and the pencil
\begin{equation}\label{eq:cp}
 \Phi_j(Y)=\nlsum_{k=1}^r K_{j,k}YK_{j,k}^*,\qquad
 \Ldm(A_1,\ldots,A_m)=\nlsum_{j=1}^m\Phi_j(A_j)
\end{equation}
satisfy $\Phi_j(I_d)=I_N$ and, for tuples $X=(A_1,\ldots,A_m)\in\HH_d(\F)^m$, the cone equalities
\[
 \ocone(\pdm,e_{d,m})=\{X:\Ldm(X)\succ0\},\qquad
 \mathcal K_{d,m}=\{X:\Ldm(X)\succeq0\}.
\]
The same fixed matrices work over $\R$ and $\C$. In addition,
\begin{equation}\label{eq:multiplier-main}
 \det\Ldm=\pdm q_{d,m},\qquad
 q_{d,m}(e_{d,m})=1,\qquad
 \mathcal K_{d,m}\subseteq\cone(q_{d,m},e_{d,m}),
\end{equation}
where $q_{d,m}$ is a homogeneous polynomial, hyperbolic with respect to $e_{d,m}$ (possibly the constant $1$).
\end{theorem}
\end{dashedbox}
Theorem~\ref{thm:main} shows that every cone $\mathcal K_{d,m}$ is spectrahedral, so the generalized
Lax conjecture holds for mixed determinants. The pencil $\Ldm$ is linear in the coordinates of
$A_1,\ldots,A_m$ alone, so these cones are spectrahedral themselves, not just spectrahedral shadows
(linear projections of spectrahedra). For complex arguments $\Ldm$ is Hermitian; realification followed by a compression gives a real symmetric pencil 
of order $2N-1$ (Remark~\ref{rmk:realification}). The factorization~\eqref{eq:multiplier-main} is the
form of spectrahedrality in Br\"and\'en's theorem on elementary symmetric
cones~\citep[(1.2)]{branden2014hyperbolicity}. Specializing to $A_j=z_jA$ answers the question raised by~\eqref{eq:rankone-intro}; we record the answer as Corollary~\ref{cor:target-intro}.
\begin{dashedbox}
\begin{corr}\label{cor:target-intro}
Let $A\in\HH_n(\F)$ with $A\succeq0$, and let $f_A(z_1,\ldots,z_n)=\eta(z_1A,\ldots,z_nA)$. Then $f_A\not\equiv0$ if and only if $A_{ii}>0$ for every $i$. Whenever the diagonal is positive, $f_A$ is hyperbolic with respect to $\one\in\R^n$, and the matrices $C_j=\Phi_j(A)\succeq0$ supplied by Theorem~\ref{thm:main} with $d=m=n$ satisfy
\begin{equation*}
 \nlsum_{j=1}^n C_j\succ0,\qquad
 \cone(f_A,\one)
 =\bigl\{z\in\R^n:\nlsum_{j=1}^n z_j C_j\succeq0\bigr\}.
\end{equation*}
The determinant of this pencil is $f_A$ times a polynomial whose hyperbolicity cone contains $\cone(f_A,\one)$. A constant congruence can normalize $\nlsum_jC_j$ to the identity.
\end{corr}
\end{dashedbox}

\subsection{Proof idea}\label{subsec:idea}
The proof is by induction on $d$. When we add one row and one column to each argument, $p_{d+1,m}$ becomes
\[
 \left(\nlsum_j\alpha_j\right)p
       -\nlsum_j D^{(j)}_{a_ja_j^*}p,
\]
where $p=\pdm(A_1,\ldots,A_m)$, and $\alpha_j$ and $a_j$ are the new diagonal entry and the new off-diagonal column of the $j$th argument, and $D^{(j)}_H$ differentiates in direction $H$ in that argument (Lemma~\ref{lem:border-identity}). We call the enlarged arguments, their polynomial, and its cone \emph{bordered}. Suppose we already have a pencil $L$ of the form~\eqref{eq:cp} that represents the cone of $p$. Since $\Phi_j(aa^*)=U_j(a)U_j(a)^*$ for a matrix $U_j(a)$ linear in $a$, and $p$ is affine along rank-one updates, we can read the derivatives at points of the open cone as extremal eigenvalues,
\[
 D^{(j)}_{aa^*}\log p
   =\lambda_{\max}\bigl(U_j(a)^*L^{-1}U_j(a)\bigr),
 \qquad U_j(a)=[K_{j,1}a\ \cdots\ K_{j,r}a].
\]
We obtain this identity by computing where a rank-one line leaves the cone, once from $p$ and once
from $L$, so it holds whatever extra factors $\det L$ carries. Dividing the bordered polynomial by
$p$, we find that the bordered cone is cut out by
$\nlsum_j\alpha_j>\nlsum_j\lambda_{\max}(U_j(a_j)^*L^{-1}U_j(a_j))$ (Lemma~\ref{lem:border-cone}). We
write the sum of these $m$ largest eigenvalues as the largest eigenvalue of one tensor sum, and a
Schur complement then turns the inequality into the next linear pencil. The determinant of this pencil
is the product of the determinants of the diagonal blocks and of the Schur complement; only the
smallest eigenvalue of the Schur complement vanishes with the bordered polynomial, so the other
factors make up the multiplier $q_{d+1,m}$. Because the step needs the factorization
$\Phi_j(aa^*)=U_j(a)U_j(a)^*$, we end by writing the new pencil again in the form~\eqref{eq:cp}, with
unital maps and new $0/1$ matrices $K_{j,k}$; the induction then continues for every $d$ and every
fixed $m$.

\subsection{Smaller pencils and further results}\label{subsec:headline-bounds}
Each step of the induction takes tensor products of the pencils of the $m$ arguments, so the order
$N_{d,m}$ of the $0/1$ pencil is doubly exponential in $d$ for fixed $m>1$
(Section~\ref{subsec:sizes}). Much smaller pencils represent the same cones, and how small depends on
which inputs the pencil must represent.

For the universal cone we adapt an argument of \citet{kummer2021spectral} and compress the pencil to
the coefficient span of a polynomial boundary kernel, a polynomial vector in the kernel of the pencil
wherever $\pdm$ vanishes; the maps stay unital and completely positive, and the order drops to
$\exp(O(md^3\log m))$ for $d,m\ge2$ (Proposition~\ref{prop:compressed-size}). Starting the induction
at order three from a pencil of order $2m^2+1$ lowers the bound further
(Propositions~\ref{prop:cubic-resolvent} and~\ref{prop:cubic-seeded-size}), at the price of
\emph{decomposable} maps for complex Hermitian arguments, that is, sums of a completely positive map
and a completely positive map composed with transposition. For the fixed-matrix cones of
Corollary~\ref{cor:target-intro} we count coefficients in the scalar variables alone, and the orders
drop further still. At $d=m=n$ the three regimes compare as follows; all entries are upper bounds.
\begin{center}
\begin{tabular}{lrr}
\toprule
Inputs represented & $\log N$ & Order at $n=4$\\
\midrule
Arbitrary matrix tuples & $O(n^4\log n)$ & $1\,639\,792$\\
Every common matrix $A$ & $O(n^3\log n)$ & $1\,639\,792$\\
One fixed matrix $A$ & $O(n^2\log n)$ & $244$\\
\bottomrule
\end{tabular}
\end{center}
At $n=4$ the first two entries coincide because the compressed bounds exceed the order of the
uncompressed pencil started at order three. For one explicit positive definite $4\times4$ matrix,
every constant compression of that pencil has order at least $244$, yet a different construction
gives order at most $20$ (supplement, Section~S2).

Compressing the universal pencil also shows that the localizations of the mixed determinant at
rank-deficient matrices have spectrahedral cones, and so do their zero-border extensions
(Proposition~\ref{prop:compression} and Corollary~\ref{cor:zero-border}). A cubic linear section of
one localization has no positive power with a definite determinantal representation, yet its cone has
a $17\times17$ pencil (Appendix~\ref{app:cubic}), so determinant multipliers are sometimes
unavoidable.

Finally, explicit polynomials of degree five answer Problems~35--38 of Borcea and Br\"and\'en
negatively, including the Schur positivity and immanant inequalities proposed there; one of them even
has a definite determinantal representation (Appendix~\ref{sec:aim}). We reproduce these examples
from our earlier investigation~\citep{Sra2026CountExMachina} because they show how far
Corollary~\ref{cor:target-intro} reaches.

\subsection{Related work}\label{sec:related-work}
Beyond the three-variable case of \citet{helton2007linear}, Kummer and Netzer prove spectrahedrality under strict hyperbolicity, meaning that $t\mapsto p(x+te)$ has only simple zeros for every $x$ outside the span of $e$~\citep[Theorem~1]{kummer2026strictly}, and OpenAI constructs a hyperbolicity cone that is not spectrahedral~\citep[Theorem~1.1 and Corollary~1.2]{openai2026nonspectrahedral}; our theorem needs no strict-hyperbolicity assumption.

The rank-one specializations in~\eqref{eq:rankone-intro} connect our construction to elementary symmetric cones. Br\"and\'en's matrix--tree construction proves their spectrahedrality~\citep[Theorem~1.2 and~(1.2)]{branden2014hyperbolicity}, and Kummer compresses that pencil to order $\binom{m}{d-1}$ for $e_d$ in $m$ variables by writing a polynomial boundary-kernel vector in a fixed coefficient space~\citep[Theorem~5.5 and Remark~5.6]{kummer2021spectral}; our size bounds come from the same argument applied to $\pdm$.

Derivative relaxations are a second comparison. Saunderson gives a pencil of order $\binom{n+1}{2}-1$ for the first derivative relaxation of the positive semidefinite cone of order $n$~\citep[Theorem~1]{saunderson2018spectrahedral}, Kummer extends spectrahedrality to all repeated derivatives of a definite determinant in one direction~\citep[Corollaries~5.12--5.13]{kummer2021spectral}, and Saunderson and Parrilo obtain size $O(\min\{k,n-k\}n^2)$ for the $k$th relaxation with auxiliary variables, as semidefinite lifts~\citep[Theorem~1]{saunderson2015polynomial}. None of them gives a pencil for independent matrix arguments.

The distinction between a polynomial and its cone runs through the paper. Br\"and\'en constructed real-zero polynomials no power of which has a definite determinantal representation~\citep[Theorem~3.3]{branden2011obstructions}, yet Kummer gave a $7\times7$ pencil for the cone of the specialized V\'amos polynomial, which has this obstruction~\citep[Section~3]{kummer2016vamos}; the sum-of-squares criterion of Kummer \etal detects such obstructions~\citep[Corollary~4.3]{kummer2015hyperbolic}, and the cubic of Appendix~\ref{app:cubic} shows the same distinction inside a limit of the mixed-determinant family. Because of such multipliers, a cone representation cannot simply be differentiated into a determinantal representation of the defining polynomial~\citep[Section~3.2]{saunderson2018spectrahedral}.

\section{The cone at order \texorpdfstring{$d+1$}{d+1} as an eigenvalue inequality}\label{sec:prelim}
Our induction passes from the cone of $\pdm$ to that of $p_{d+1,m}$, and we need two facts: deleting a common index from every argument interlaces the zeros along the reference direction, and the signs of those zeros determine the cone. We then add one row and one column to every argument, which reverses the deletion, and reduce membership in the cone at order $d+1$ to a scalar inequality whose terms are largest eigenvalues of compressions of $L^{-1}$, where $L$ is the pencil at order $d$.

\subsection{Interlacing and the positive component}\label{subsec:hyperbolicity}
The first fact comes from the stability theory of \citet{branden08}.

\begin{lemma}[Real-rootedness and common-index deletion]\label{lem:interlacing}
For all $A_1,\ldots,A_m\in\HH_d(\F)$, the polynomial $t\mapsto\eta(A_1+tI_d,A_2,\ldots,A_m)$ is monic of degree $d$ and real-rooted. If $i\in[d]$, the monic polynomial of degree $d-1$ obtained by deleting row and column $i$ from every argument interlaces it.
\end{lemma}
\begin{proof}
We apply \citet[Theorem~2.6]{branden08} to the pencils $A_1+tI_d,A_2,\ldots,A_m$. The same theorem gives stability of $f(t)+v g(t)$, where $f(t)=\eta(A_1+tI_d,A_2,\ldots,A_m)$ and $g$ is the common-index deletion; we substitute $v=\sqrt{-1}$ and apply the Hermite--Biehler theorem~\citep[Theorem~2.8]{branden08} to obtain interlacing. Both leading coefficients are $1$, from the partitions that assign every remaining index to the first argument, so neither polynomial is zero; for $d=1$ the deletion is the constant $1$.
\end{proof}

Lemma~\ref{lem:interlacing} holds equally with any one argument shifted by $tI_d$, so each tuple with $I_d$ in exactly one position and zeros elsewhere is a direction of hyperbolicity. These directions lie in the same positive component, since $\eta$ is positive on the segments joining them, by the identity
\begin{equation}\label{eq:scalar-identities}
 \eta(t_1I_d,\ldots,t_mI_d)=(t_1+\cdots+t_m)^d.
\end{equation}

The second fact lets us locate the positive component by the zeros along the reference direction.

\begin{lemma}[Root description of the component]\label{lem:component}
Let $p$ be homogeneous and hyperbolic with respect to $e$, with $p(e)>0$. The connected component of $\{p>0\}$ containing $e$ is
\[
 \ocone(p,e)=\{x:\text{all zeros of }p(x+te)\text{ are strictly negative}\}.
\]
Replacing ``strictly negative'' by ``nonpositive'' gives its closure, whose boundary lies in $\{p=0\}$.
\end{lemma}
\begin{proof}
The leading coefficient of $p(x+te)$ is the fixed nonzero number $p(e)$, so the real zeros vary continuously as a multiset. Along a path where $p(x)\ne0$, no zero crosses $0$, so the number of negative zeros is locally constant. Conversely, if all zeros are negative, we move along the segment $(1-u)x+ue$, $0\le u\le1$, whose zeros are $(1-u)\rho_i-u$, where $\rho_i$ are the zeros at $x$; they remain negative. These two observations prove the component assertion. For the closure we add $\eps e$, which shifts every zero by $-\eps$. At a boundary point at least one zero is $0$, proving the last assertion.
\end{proof}

In the induction, the convexity of each cone comes from its representing pencil. Lemmas~\ref{lem:interlacing} and~\ref{lem:component} hold with repeated zeros and with a nontrivial lineality space.

\subsection{Exit thresholds as extremal eigenvalues}\label{sec:border}
We now write the bordered polynomial in terms of directional derivatives of the mixed determinant, and show that the mixed determinant is affine along rank-one updates. For a polynomial in matrix arguments, we write $D_H^{(j)}$ for differentiation in direction $H$ in argument $j$, on the underlying real vector space, also in the Hermitian case.

\begin{lemma}[Polynomial identities]\label{lem:border-identity}
Let $p=\pdm(A_1,\ldots,A_m)$ and write
\[
 A_j^+=\begin{pmatrix}A_j&a_j\\a_j^*&\alpha_j\end{pmatrix},\qquad
 \sigma=\nlsum_{j=1}^m\alpha_j.
\]
The bordered polynomial then satisfies
\begin{equation}\label{eq:border}
 p_{d+1,m}(A_1^+,\ldots,A_m^+)
     =\sigma p-\nlsum_{j=1}^mD_{a_ja_j^*}^{(j)}p.
\end{equation}
For every $j$ and every $a\in\F^d$,
\begin{equation}\label{eq:rankone-affinity}
 p(A_1,\ldots,A_j-taa^*,\ldots,A_m)
       =p-tD_{aa^*}^{(j)}p.
\end{equation}
Both identities hold for all Hermitian arguments.
\end{lemma}
\begin{proof}
By the matrix determinant lemma, every principal determinant is affine under a rank-one update, also at singular matrices; each term of the mixed determinant is therefore affine in $t$, which proves~\eqref{eq:rankone-affinity}. For~\eqref{eq:border} we look at a partition contributing to the bordered polynomial: the new index lies in exactly one part, say the $j$th, and expanding its bordered determinant gives the old determinant times $\alpha_j$, minus the directional derivative in direction $a_j[S_j]a_j[S_j]^*$, which is zero for an empty part. Summing over the old partitions and over $j$, we get~\eqref{eq:border}.
\end{proof}

Along a rank-one update $x-th$ the polynomial is affine in $t$, say $p(x-th)=p(x)(1-t\theta)$, so the line leaves the open cone at most once, at $t=1/\theta$; we call $\theta$ the \emph{exit threshold}. To compute $\theta$ we need only a pencil that represents the cone, so we work with any hyperbolic polynomial.

\begin{dashedbox}[supporting]
\begin{lemma}[Rank-one exit threshold]\label{lem:exit}
Let $p$ be homogeneous and hyperbolic with respect to $e$, with $p(e)>0$, and suppose that a homogeneous linear pencil $L$ satisfies $\ocone(p,e)=\{L\succ0\}$ and $\cone(p,e)=\{L\succeq0\}$. Fix $x\in\ocone(p,e)$ and a direction $h$ such that $L(h)=UU^*$ for a rectangular matrix $U$ and $p(x-th)=p(x)(1-t\theta)$ for all real $t$. Then
\begin{equation*}
 \theta=\lambda_{\max}(U^*L(x)^{-1}U)
        =\norm{L(x)^{-1/2}U}_{\op}^2\ge0.
\end{equation*}
The value is understood to be zero when $U=0$.
\end{lemma}
\end{dashedbox}
\begin{proof}
Put $\gamma=\lambda_{\max}(U^*L(x)^{-1}U)$. If $\gamma>0$, the matrix $L(x)-tUU^*$ is positive definite for $t<1/\gamma$, and positive semidefinite and singular at $t=1/\gamma$. By the cone equalities, the latter point lies on the boundary of $\cone(p,e)$, so Lemma~\ref{lem:component} gives $p(x-h/\gamma)=0$, and we conclude that $\theta=\gamma$. If $\gamma=0$, then $U=0$ and $L(x-th)\succ0$ for every real $t$; the affine function $p(x)(1-t\theta)$ is then positive for every real $t$, which forces $\theta=0$.
\end{proof}

The threshold is an extremal eigenvalue, not the trace $\tr(L^{-1}L(h))$, which is the derivative of $\log\det L$ and so also picks up the extra factors of $\det L$. To apply Lemma~\ref{lem:exit} to the mixed determinant, we assume that Theorem~\ref{thm:main}, except for the multiplier assertion, holds at matrix order $d$. At $X=(A_1,\ldots,A_m)$ in the open cone $\mathcal K_{d,m}^{\circ}=\ocone(\pdm,e_{d,m})$ we put $L=\Ldm(X)$ and introduce $U_j(a)=[K_{j,1}a\ \cdots\ K_{j,r}a]$ and $\theta_j(X;a)=D_{aa^*}^{(j)}p(X)/p(X)$. Then $\Phi_j(aa^*)=U_j(a)U_j(a)^*$, and Lemma~\ref{lem:exit} with~\eqref{eq:rankone-affinity} yields
\begin{equation}\label{eq:theta}
 \theta_j(X;a)
 =\lambda_{\max}\bigl(U_j(a)^*L^{-1}U_j(a)\bigr)\ge0.
\end{equation}

By~\eqref{eq:border}, the bordered polynomial equals $p(X)\bigl(\sigma-\nlsum_j\theta_j(X;a_j)\bigr)$ at such $X$, and we now show that its positivity, together with $X\in \mathcal K_{d,m}^{\circ}$, characterizes the bordered open cone.

\begin{lemma}[The bordered positive component]\label{lem:border-cone}
Under the preceding induction hypothesis,
\begin{equation*}
 (A_1^+,\ldots,A_m^+)\in \mathcal K_{d+1,m}^{\circ}
 \quad\Longleftrightarrow\quad
 \begin{cases}
 X\in \mathcal K_{d,m}^{\circ},\\
 \displaystyle\sigma>\nlsum_{j=1}^m\theta_j(X;a_j).
 \end{cases}
\end{equation*}
\end{lemma}
\begin{proof}
For necessity, we shift the first bordered argument by $tI_{d+1}$; the bordered polynomial then has all zeros strictly negative. Deleting the new index and applying Lemma~\ref{lem:interlacing}, we find that the corresponding zeros for the base tuple are also strictly negative, so $X\in \mathcal K_{d,m}^{\circ}$. Since $p(X)>0$, identity~\eqref{eq:border} and positivity of the bordered polynomial give the displayed inequality.

For sufficiency, we write $\Theta=\nlsum_j\theta_j(X;a_j)\ge0$ and scale all border vectors by $u\in[0,1]$, which changes the bordered polynomial to $p(X)(\sigma-u^2\Theta)>0$; so the tuple is connected through positive values to one with all border vectors zero, where the polynomial is $\sigma p(X)$ with $\sigma>0$. We then move $X$ to $e_{d,m}$ along the straight segment inside $\{\Ldm\succ0\}$, and finally move $(\alpha_1,\ldots,\alpha_m)$ to $(1,0,\ldots,0)$ along a segment whose sum remains positive. The path we have traced ends at $e_{d+1,m}$ and never crosses a zero, so the bordered tuple lies in the component of $e_{d+1,m}$, that is, in $\mathcal K_{d+1,m}^{\circ}$.
\end{proof}

\section{Proof of the main theorem}\label{sec:induction}
We now construct the pencils of Theorem~\ref{thm:main} by induction on $d$, with $m\ge1$ fixed. We saw in Lemma~\ref{lem:border-cone} that the inequality $\sigma>\nlsum_j\theta_j(X;a_j)$ cuts out the bordered cone over $\mathcal K_{d,m}^{\circ}$, and in~\eqref{eq:theta} that each $\theta_j$ is the largest eigenvalue of a compression of $L^{-1}$. Since the eigenvalues of a tensor sum add, we can express this inequality linearly by a block pencil whose Schur complement is $\sigma$ minus a tensor sum of these compressions.

For $d=1$ we take $N=r=1$ and $K_{j,1}=[1]$; then $p_{1,m}(A_1,\ldots,A_m)=A_1+\cdots+A_m$ and $L_{1,m}=[A_1+\cdots+A_m]$ for real scalars $A_j$, and the cone equalities and the Kraus form~\eqref{eq:cp} hold.

Suppose these assertions hold at order $d$, with sizes $N,r$. We set
\[
 E=\F^r,\qquad V=\F^N,\qquad T=E^{\otimes m},
 \qquad V_j=E^{\otimes(j-1)}\otimes V\otimes E^{\otimes(m-j)},
\]
where a zero tensor power denotes the scalar field, and we define the operators $\mathcal D_j$ on $V_j$ and $\mathcal U_j$ from $T$ to $V_j$ by
\begin{align*}
 \mathcal D_j&=I_r^{\otimes(j-1)}\otimes L\otimes I_r^{\otimes(m-j)},\\
 \mathcal U_j&=I_r^{\otimes(j-1)}\otimes U_j(a_j)\otimes I_r^{\otimes(m-j)}.
\end{align*}
On $V_1\oplus\cdots\oplus V_m\oplus T$ we put
\begin{equation}\label{eq:new-pencil}
 L^+=
 \begin{pmatrix}
 \mathcal D_1&0&\cdots&0&\mathcal U_1\\
 0&\mathcal D_2&\cdots&0&\mathcal U_2\\
 \vdots&\vdots&\ddots&\vdots&\vdots\\
 0&0&\cdots&\mathcal D_m&\mathcal U_m\\
 \mathcal U_1^*&\mathcal U_2^*&\cdots&\mathcal U_m^*&\sigma I_{r^m}
 \end{pmatrix}.
\end{equation}
The entries are linear in the bordered inputs, since the tensor factors are constant identity matrices.

\subsection{The new pencil cuts out the open cone}
If $L\succ0$, the Schur complement of the block $\operatorname{diag}(\mathcal D_1,\ldots,\mathcal D_m)$ in~\eqref{eq:new-pencil} is
\begin{equation}\label{eq:schur}
 \sigma I_{r^m}
 -\nlsum_{j=1}^m I_r^{\otimes(j-1)}\otimes H_j\otimes I_r^{\otimes(m-j)},
 \qquad H_j=U_j(a_j)^*L^{-1}U_j(a_j).
\end{equation}
For Hermitian matrices $H_1,\ldots,H_m$ we have, by the spectral theorem,
\begin{equation}\label{eq:tensorsum}
 \lambda_{\max}\left(
 \nlsum_{j=1}^m I^{\otimes(j-1)}\otimes H_j\otimes I^{\otimes(m-j)}\right)
   =\nlsum_{j=1}^m\lambda_{\max}(H_j),
\end{equation}
because tensor products of eigenvectors form an eigenbasis and their eigenvalues add; the embeddings of the $H_j$ act on different tensor factors and commute.

If $L^+\succ0$, each principal block $\mathcal D_j$ is positive definite, so $L\succ0$ and the Schur complement~\eqref{eq:schur} is defined. Combining~\eqref{eq:theta}--\eqref{eq:tensorsum} with Lemma~\ref{lem:border-cone}, we obtain
\begin{equation}\label{eq:open-induction}
 L^+\succ0\quad\Longleftrightarrow\quad
 (A_1^+,\ldots,A_m^+)\in \mathcal K_{d+1,m}^{\circ}.
\end{equation}

\subsection{The new maps are again unital and completely positive}\label{subsec:new-maps}
For the induction to continue, we must write $L^+$ again as a sum of unital completely positive maps in Kraus form, one for each bordered argument. Let $e_1,\ldots,e_r$ be the standard basis of $E$. For a multi-index
$\mathbf k=(k_1,\ldots,k_m)\in[r]^m$, write
$e_{\mathbf k}=e_{k_1}\otimes\cdots\otimes e_{k_m}$.
For each argument $c\in[m]$, each top block $j\in[m]$, and each $\mathbf k$, we define $K^+_{c,j,\mathbf k}\colon\F^{d+1}\to V_1\oplus\cdots\oplus V_m\oplus T$ by specifying its only possibly nonzero components on $(v,s)\in\F^d\oplus\F$:
\begin{align}
 (K^+_{c,j,\mathbf k}(v,s))_{V_j}
   &=e_{k_1}\otimes\cdots\otimes K_{c,k_j}v\otimes\cdots\otimes e_{k_m},\label{eq:newK-top}\\
 (K^+_{c,j,\mathbf k}(v,s))_T
   &=\delta_{cj}\,s\,e_{\mathbf k}.\label{eq:newK-bottom}
\end{align}
In~\eqref{eq:newK-top}, the old Kraus operator $K_{c,k_j}$ occupies tensor position $j$. We define
\begin{equation}\label{eq:newPhi}
 \Phi_c^+(Y)=\nlsum_{j=1}^m\ \sum_{\mathbf k\in[r]^m}
      K^+_{c,j,\mathbf k}Y (K^+_{c,j,\mathbf k})^*
\end{equation}
and check every block of $\Phi_c^+(A_c^+)$. The $j$th top diagonal block is $I_r^{\otimes(j-1)}\otimes\Phi_c(A_c)\otimes I_r^{\otimes(m-j)}$, and the bottom diagonal block is $\alpha_cI_{r^m}$, since only $j=c$ contributes there. The only nonzero top-to-bottom block occurs for $j=c$ and equals $I_r^{\otimes(c-1)}\otimes[K_{c,1}a_c\ \cdots\ K_{c,r}a_c]\otimes I_r^{\otimes(m-c)}=\mathcal U_c$. Each $K^+_{c,j,\mathbf k}$ maps into one top block, so all couplings between distinct top blocks vanish. Summing over $c$, we obtain
\begin{equation}\label{eq:cp-closes}
 L^+(A_1^+,\ldots,A_m^+)=\nlsum_{c=1}^m\Phi_c^+(A_c^+).
\end{equation}
When we substitute $Y=I_{d+1}$ in~\eqref{eq:newPhi}, every top diagonal block and the bottom diagonal block become the identity, while all off-diagonal blocks vanish, so
\begin{equation}\label{eq:unital-closes}
 \Phi_c^+(I_{d+1})=I\qquad(1\le c\le m).
\end{equation}
Every new Kraus operator has entries in $\{0,1\}$, because the old ones do.

\noindent\textbf{Closed cones and completion of the induction.} By~\eqref{eq:unital-closes}, $L^+(e_{d+1,m})=I$, so for every tuple $X^+$ of order $d+1$, $L^+(X^+)\succeq0$ implies $L^+(X^++\eps e_{d+1,m})\succ0$ for every $\eps>0$, and taking closures in~\eqref{eq:open-induction} we obtain the closed-cone equality in both directions. The new maps~\eqref{eq:newPhi} satisfy exactly the induction hypothesis at order $d+1$, and we have proved the Kraus form~\eqref{eq:cp} and the cone equalities for every $d,m$, including $m=1$. It remains to prove the multiplier assertion~\eqref{eq:multiplier-main}.

\subsection{Irreducibility and the determinant multiplier}\label{sec:multiplier}
For the factorization~\eqref{eq:multiplier-main} we show that the bordered polynomial divides $\det L^+$. As a polynomial in $\sigma$, the bordered polynomial is linear, and $\det L^+$ vanishes at its root, so the division is exact over the field of fractions of the other coordinates; it is exact in the polynomial ring because the bordered polynomial is primitive in $\sigma$, which we derive from the irreducibility of $p_{d,m}$. All polynomial rings in this subsection are over $\R$, in the independent real coordinates of the symmetric or Hermitian inputs.

\begin{prop}\label{prop:irreducible}
The universal polynomial $p_{d,m}$ is irreducible for every $d,m\ge1$.
\end{prop}
\begin{proof}
For $d=1$ the polynomial is linear. Assume $p=p_{d,m}$ irreducible, and write the next polynomial as $p_{d+1,m}=\sigma p-R$ with $R=\nlsum_jD_{a_ja_j^*}^{(j)}p$, using $\sigma$ and $m-1$ independent differences of the $\alpha_j$ as coordinates in place of $(\alpha_1,\ldots,\alpha_m)$; this invertible linear change does not affect irreducibility. We view $R$ as a polynomial in the real border coordinates; each coefficient is homogeneous of degree $d-1$ in the base matrix coordinates, whereas $p$ has degree $d$, and at least one coefficient is nonzero, as we see by taking $A_1=I_d$, the other base arguments zero, $a_1\ne0$, and the other borders zero. Since $p$ is prime in the base polynomial ring and remains prime after adjoining variables, $\gcd(p,R)=1$, so $\sigma p-R$ is primitive as a polynomial in $\sigma$. Being linear in $\sigma$, it is irreducible over the fraction field of the other coordinates, and Gauss's lemma gives irreducibility in the polynomial ring.
\end{proof}

\begin{prop}\label{prop:multiplier}
The pencil constructed in this section satisfies~\eqref{eq:multiplier-main}.
\end{prop}
\begin{proof}
The assertion is immediate at $d=1$. At the next step we use the notation of the preceding proof. For every interior base tuple and every choice of border vectors, setting $\sigma=R/p$ makes the Schur complement~\eqref{eq:schur} positive semidefinite and singular, by~\eqref{eq:theta} and~\eqref{eq:tensorsum}, so $\det L^+$ vanishes at this value of $\sigma$. The vanishing holds on an open set of the other real coordinates, so it is a rational-function identity; therefore $\sigma p-R$ divides $\det L^+$ over the fraction field, and primitivity, as in Proposition~\ref{prop:irreducible}, gives divisibility in the polynomial ring.

We write $\det L_{d,m}=p_{d,m}q_{d,m}$. The quotient is homogeneous of degree $N_{d,m}-d$, and evaluating at $e_{d,m}$ we get $q_{d,m}(e_{d,m})=1$. On every real line parallel to $e_{d,m}$, the determinant is the characteristic polynomial, up to sign convention, of a Hermitian matrix shifted by $tI$, so all zeros of the quotient are real, and the quotient is hyperbolic, with the usual convention for a nonzero constant. On $\mathcal K_{d,m}^{\circ}$ both $p_{d,m}$ and the determinant are positive, so $q_{d,m}>0$ there, and connectedness and closure give the cone inclusion in~\eqref{eq:multiplier-main}.
\end{proof}

With Proposition~\ref{prop:multiplier} the proof of Theorem~\ref{thm:main} is complete.

\subsection{Doubly exponential orders, sparse Kraus operators}\label{subsec:sizes}
The sizes satisfy
\begin{align*}
 N_{1,m}=r_{1,m}&=1,\\
 N_{d+1,m}&=mN_{d,m}r_{d,m}^{m-1}+r_{d,m}^m,\\
 r_{d+1,m}&=m r_{d,m}^m.
\end{align*}
For $m>1$ these recursions solve to
\begin{equation*}
 r_{d,m}=m^{(m^{d-1}-1)/(m-1)},\qquad
 N_{d,m}=\frac{m+d-1}{m}\,r_{d,m}.
\end{equation*}
For $m=1$, $r_{d,1}=1$ and $N_{d,1}=d$. The first few values of the $0/1$ construction are as follows; in Section~\ref{sec:smaller} we represent the same cones by much smaller pencils.

\begin{center}
\begin{tabular}{rrrr}
\toprule
$d$&$m$&$N_{d,m}$&$r_{d,m}$\\
\midrule
2&2&3&2\\
2&3&4&3\\
3&2&16&8\\
3&3&135&81\\
3&4&1536&1024\\
4&2&320&128\\
\bottomrule
\end{tabular}
\end{center}

The Kraus operators are sparse in a strong sense: each $K_{j,k}$ is a partial coordinate embedding, that is, each column is either zero or a standard basis vector and the nonzero columns are distinct, and for a fixed $j$ the row supports of these maps partition all output coordinates. We check the partition property by induction. It holds in the base case; in~\eqref{eq:newK-top}, a top output coordinate determines the tensor positions other than $j$ and, by the old partition property, its unique old map and input coordinate; and in~\eqref{eq:newK-bottom}, a bottom coordinate occurs exactly once, for $j=c$. Consequently, each $\Phi_j$ is, after a fixed permutation of output coordinates, a direct sum of principal-submatrix maps, which gives us a purely combinatorial check of unitality.

As a check on the bookkeeping, we verified the identity~\eqref{eq:cp-closes} and the unitality~\eqref{eq:unital-closes} in exact arithmetic for $(d,m)=(1,1),(3,1),(3,2),(3,3),(3,4),(4,2)$, by comparing the maps~\eqref{eq:newK-top}--\eqref{eq:newK-bottom} with a separate assembly of the block formula~\eqref{eq:new-pencil}.

\section{Smaller pencils for the universal cone}\label{sec:smaller}
The tensor products of Section~\ref{sec:induction} make $N_{d,m}$ doubly exponential in $d$ for fixed $m>1$. To do better, we take a polynomial boundary kernel of the pencil and compress the pencil to the coefficient span of that kernel: enough of the kernel survives the compression to preserve the cone, so we obtain a smaller pencil for the same cone. We adapt this argument from \citet[Theorem~5.5]{kummer2021spectral}, who used it for elementary symmetric cones. Compression preserves unital complete positivity when the original maps have it, but the compressed Kraus operators need no longer have entries in $\{0,1\}$.

\subsection{Compression and the general bound}
We first show that a compression which retains an open part of the boundary cannot enlarge the cone, because the defining polynomial is irreducible.

\begin{lemma}[Compression by a boundary kernel]\label{lem:kernel-compression}
Let $p$ be an irreducible homogeneous polynomial hyperbolic with respect to $e$, with $p(e)=1$. Suppose that a homogeneous Hermitian pencil $L$ satisfies $L(e)=I$, $\ocone(p,e)=\{L\succ0\}$, and $\cone(p,e)=\{L\succeq0\}$. Let $v$ be a \emph{polynomial boundary kernel} of $L$, that is, a homogeneous polynomial vector such that $L(y)v(y)=0$ whenever $p(y)=0$, and suppose that $v$ is nonzero on a relatively open smooth part of the cone boundary, a \emph{boundary chart}. If $Q$ has orthonormal columns whose span contains the coefficient vectors of $v$, then
\[
 \ocone(p,e)=\{Q^*LQ\succ0\},\qquad
 \cone(p,e)=\{Q^*LQ\succeq0\},\qquad
 p\mid\det(Q^*LQ).
\]
Compression also preserves the unital complete positivity of each map in a sum $L=\nlsum_j\Phi_j$.
\end{lemma}
\begin{proof}
Put $M=Q^*LQ$. At a boundary point where $v(y)\ne0$, we write $v(y)=Qu$; then $u\ne0$ and $M(y)u=0$. Thus $\det M$ vanishes on a relatively open smooth subset of $\{p=0\}$, and irreducibility gives $p\mid\det M$. The convex set $\{M\succ0\}$ contains $e$ and avoids $\{p=0\}$, so it lies in $\ocone(p,e)$; compression from $L\succ0$ gives the reverse inclusion. Since $M(e)=I$, adding $\eps e$ gives equality of the closures. Finally, compressing a Kraus representation replaces $K_{j,k}$ by $Q^*K_{j,k}$, and $Q^*\Phi_j(I)Q=I$.
\end{proof}

We bound the dimension of the coefficient span by the degree of the polynomial vector.

\begin{dashedbox}[supporting]
\begin{prop}[Compressed size]\label{prop:compressed-size}
For $d\ge2$ and $m\ge1$, set
\begin{equation}\label{eq:kernel-degree}
 \begin{split}
 n_{d,m}&=d+md(d-1),\qquad
 \kappa_{d,m}=\frac{m(2m)^{d-2}+m-1}{2m-1},\\
 \delta_{d,m}&=(d-1)\kappa_{d,m}.
 \end{split}
\end{equation}
There is a pencil $\nlsum_j\widehat\Phi_j(A_j)$ with real algebraic Kraus operators and unital completely positive maps that satisfies the cone equalities and the multiplier assertion of Theorem~\ref{thm:main}, of order
\begin{equation}\label{eq:compressed-size}
 \widehat N_{d,m}\le
 \min\left\{N_{d,m},
       \binom{n_{d,m}+\delta_{d,m}-2}{n_{d,m}-2}\right\}.
\end{equation}
The same Kraus operators work over $\R$ and $\C$, and each map needs at most $d\widehat N_{d,m}$ Kraus operators. In particular, $\log\widehat N_{d,m}=O(md^3\log m)$ for $d,m\ge2$. For real symmetric arguments alone, $n_{d,m}$ in~\eqref{eq:compressed-size} can be replaced by $d+md(d-1)/2$.
\end{prop}
\end{dashedbox}
The binomial counts the monomials of a polynomial boundary kernel of
degree $\delta_{d,m}$ that omits the last diagonal sum; when we clear
denominators in the induction, it is the rank-one degree $\kappa_{d,m}$
of the kernel, not its full degree, that enters. We construct this kernel
and prove Proposition~\ref{prop:compressed-size} in
Appendix~\ref{app:kernel-construction}. The bound is increasing in the
kernel degree, so starting the induction from a kernel of lower degree
improves every later bound. At order three, grouping the kernel
coordinates already gives a smaller completely positive bound
(Appendix~\ref{subsec:order-three-size}), and once we allow transposed
summands we find a pencil of order $2m^2+1$ with a kernel of degree two.
We construct it next and continue the induction from it in
Section~\ref{subsec:cubic-seeded-size}.

\subsection{A quadratic-size pencil at order three}\label{subsec:cubic-resolvent}
At order three we can write the logarithmic derivatives of $p_{2,m}$
along rank-one updates as quadratic forms in the inverses of pencils of
order $2m$, one for each argument, and then a single scalar Schur
complement expresses their sum. The price we pay is that for Hermitian
arguments the construction uses maps of the form
\begin{equation}\label{eq:decomposable-map}
 \Psi(Y)=\nlsum_k K_kYK_k^*+\nlsum_\ell J_\ell Y^{\mathsf T}J_\ell^*.
\end{equation}
The first sum is completely positive, and the second is a completely
positive map composed with transposition; we call such a sum
\emph{decomposable}. On real symmetric inputs, replacing
$Y^{\mathsf T}$ by $Y$ gives a completely positive extension of the same
map.

\begin{prop}[Cubic resolvent construction]\label{prop:cubic-resolvent}
For $m\ge2$, the universal cone $\mathcal K_{3,m}$ has a pencil
$L^{\mathrm{res}}_{3,m}=\nlsum_j\Psi_j(A_j)$ of order $2m^2+1$, with
$\Psi_j(I_3)=I$. Each map has a representation~\eqref{eq:decomposable-map}
using $\rho_{3,m}=2m^2-3m+2$ matrices in total, all with entries in $\{0,1,-1\}$. The same matrices
work over both fields. The Hermitian coefficient maps are decomposable;
on real symmetric inputs they admit unital completely positive extensions.
The pencil satisfies the open- and closed-cone equalities and the
multiplier assertion of Theorem~\ref{thm:main}. It also has a polynomial boundary kernel of degree two
and of degree at most one along rank-one updates in one argument.
\end{prop}
\begin{proof}
We write the arguments as
\[
 \begin{gathered}
 A_j^+=\begin{pmatrix}A_j&a_j\\a_j^*&\alpha_j\end{pmatrix},\\
 s=\nlsum_j(A_j)_{11},\quad t=\nlsum_j(A_j)_{22},\quad
 b_j=(A_j)_{12},\quad \sigma=\nlsum_j\alpha_j,
 \end{gathered}
\]
and put $p=st-\nlsum_j|b_j|^2$ and
\[
 q_j=a_j^*\begin{pmatrix}t&-b_j\\-\overline b_j&s\end{pmatrix}a_j,
\]
so that Lemma~\ref{lem:border-identity} gives
$p_{3,m}(A_1^+,\ldots,A_m^+)=\sigma p-\nlsum_jq_j$.

To realize $q_j/p$, we let $\beta_j=(b_k)_{k\ne j}$ be a column vector and set
\begin{equation}\label{eq:cubic-resolvent}
 B_j=\begin{pmatrix}
 b_j&\beta_j^{\mathsf T}\\
 -\overline{\beta_j}&\overline b_jI_{m-1}
 \end{pmatrix},\qquad
 C_j=\begin{pmatrix}sI_m&B_j\\B_j^*&tI_m\end{pmatrix}.
\end{equation}
We index each half of $\F^{2m}$ by $0$ and the indices $k\ne j$, write
its coordinate vectors as $e_0^+,e_k^+,e_0^-,e_k^-$, and define
\[
 K_j=[e_0^+\ e_0^-],\qquad
 P_j=\begin{pmatrix}
 t&-b_j\\
 0&\overline{\beta_j}\\
 -\overline b_j&s\\
 -\overline{\beta_j}&0
 \end{pmatrix}.
\]
Direct multiplication gives
\begin{equation}\label{eq:cubic-resolvent-identities}
 \begin{gathered}
 C_jP_j=pK_j,\qquad
 K_j^*P_j=\begin{pmatrix}t&-b_j\\-\overline b_j&s\end{pmatrix},\\
 B_jB_j^*=B_j^*B_j=
 \begin{pmatrix}
 \nlsum_k|b_k|^2&0\\
 0&|b_j|^2I_{m-1}+\overline{\beta_j}\beta_j^{\mathsf T}
 \end{pmatrix}.
 \end{gathered}
\end{equation}
Thus $\norm{B_j}_{\op}^2=\nlsum_k|b_k|^2$, so $C_j\succ0$ exactly when
$s>0$ and $p>0$, and on this set
$a_j^*K_j^*C_j^{-1}K_ja_j=q_j/p$.

We now assemble the cubic pencil
\begin{equation}\label{eq:cubic-resolvent-pencil}
 L^{\mathrm{res}}_{3,m}=
 \begin{pmatrix}
 C_1&&&K_1a_1\\
 &\ddots&&\vdots\\
 &&C_m&K_ma_m\\
 a_1^*K_1^*&\cdots&a_m^*K_m^*&\sigma
 \end{pmatrix}.
\end{equation}
Its Schur complement is $\sigma-\nlsum_jq_j/p$, so by
Lemma~\ref{lem:border-cone} its positive definite locus is
$\mathcal K_{3,m}^{\circ}$, and since $L^{\mathrm{res}}_{3,m}(e_{3,m})=I$, adding
$\varepsilon e_{3,m}$ gives the closed-cone equality. From the last
identity in~\eqref{eq:cubic-resolvent-identities} we also get
$\det C_j=p^2(st-|b_j|^2)^{m-2}$, hence
\begin{equation}\label{eq:cubic-resolvent-determinant}
 \det L^{\mathrm{res}}_{3,m}
 =p_{3,m}\,p^{2m-1}\nlprod_{j=1}^m(st-|b_j|^2)^{m-2}.
\end{equation}
Each quadratic factor is hyperbolic in direction $e_{3,m}$ and positive
on the open cone, so the multiplier has the required cone containment.

Next we read off explicit maps of the form~\eqref{eq:decomposable-map}
from the coordinate pairs in~\eqref{eq:cubic-resolvent}; we refer to the
argument index of a summand as its \emph{color}. For color $j$ in block
$C_j$ we use the embedding $[e_0^+\ e_0^-]$ on $Y$ and the embeddings
$[e_k^+\ e_k^-]$, $k\ne j$, on $Y^{\mathsf T}$, and for color $c\ne j$ we
use $[e_0^+\ e_c^-]$ on $Y$ and $[e_c^+\ {-e_0^-}]$ on $Y^{\mathsf T}$. The remaining coordinates receive
$Y_{11}e_k^+(e_k^+)^*$ and $Y_{22}e_k^-(e_k^-)^*$ for $k\ne j,c$,
whose Kraus operators have one nonzero column. Each color
map sends $I_2$ to $I_{2m}$.

In the direct sum of the $C_j$ we extend the embedding
$[e_0^+\ e_0^-]$ in block $C_c$ by sending the third input coordinate
to the final coordinate of~\eqref{eq:cubic-resolvent-pencil}, producing
the terms $K_ca_c$ and $\alpha_c$, and give every other embedding a zero
third column. The resulting color-$c$ map uses
$m+(m-1)(2m-2)=\rho_{3,m}$ signed matrices, whose row supports are
disjoint, so the map is unital.
For complex inputs these maps are not completely positive: compressing
to a transposed pair in $C_c$ and restricting to the first two input
coordinates gives the transposition of $2\times2$ matrices.

From the resolvent identities we obtain the polynomial boundary kernel
\begin{equation}\label{eq:cubic-resolvent-kernel}
 v^{\mathrm{res}}_3
   =\begin{pmatrix}-P_1a_1\\\vdots\\-P_ma_m\\p\end{pmatrix},\qquad
 L^{\mathrm{res}}_{3,m}v^{\mathrm{res}}_3=p_{3,m}e_{\mathrm{last}},
\end{equation}
where $e_{\mathrm{last}}$ is the last coordinate vector. This kernel has
degree two, omits $\sigma$, and is nonzero on the smooth boundary chart
with interior base, because its last coordinate is $p>0$.

To check its rank-one degree, we update the bordered argument $A_c^+$ in
direction $(x,y,\gamma)(x,y,\gamma)^*$: the base changes linearly, and
among the border vectors only $a_c$ changes, to $a_c+\tau(x,y)^{\mathsf T}\overline\gamma$.
The possible $\tau^2$ coefficient vanishes because $P_c\bigl(h_c(x,y)\bigr)(x,y)^{\mathsf T}=0$, where $h_c(x,y)$ is the base rank-one direction; indeed, the only nonzero rows of $P_c\bigl(h_c(x,y)\bigr)$
form the matrix
$\left(\begin{smallmatrix}|y|^2&-x\overline y\\
-\overline x y&|x|^2\end{smallmatrix}\right)$, which annihilates
$(x,y)^{\mathsf T}$. The other blocks are linear in $\tau$, and $p$ is
rank-one affine, so $v^{\mathrm{res}}_3$ has rank-one degree at most
one, while $e_{\mathrm{last}}$ is constant.
\end{proof}

\subsection{Starting the induction at order three}\label{subsec:cubic-seeded-size}
We now continue the induction of Section~\ref{sec:induction} from the
cubic pencil instead of from order one, so that its smaller kernel
degree carries over to every later order. For Hermitian arguments we keep track of
the decomposition into completely positive and transposed completely
positive summands.

\begin{dashedbox}[supporting]
\begin{prop}[Continuation from the cubic pencil]\label{prop:cubic-seeded-size}
For $d\ge3$ and $m\ge2$, define
\[
 \check\kappa_{d,m}
   =\frac{m(2m)^{d-3}+m-1}{2m-1},\qquad
 \check\delta_{d,m}=(d-1)\check\kappa_{d,m}.
\]
Put $\rho_{3,m}=2m^2-3m+2$ and
\[
 \rho_{d,m}
   =m^{(m^{d-3}-1)/(m-1)}\rho_{3,m}^{\,m^{d-3}},\qquad
 S_{d,m}=\rho_{d,m}
      \left(\frac{2m^2+1}{\rho_{3,m}}+\frac{d-3}{m}\right).
\]
There is a homogeneous pencil $\check L_{d,m}=\nlsum_j\check\Phi_j(A_j)$
representing the open and closed cones of
$p_{d,m}$ and satisfying the multiplier assertion of
Theorem~\ref{thm:main}, whose order obeys
\begin{equation}\label{eq:cubic-seeded-size}
 \check N_{d,m}\le
 \min\left\{\widehat N_{d,m},S_{d,m},
       \binom{n_{d,m}+\check\delta_{d,m}-2}{n_{d,m}-2}\right\}.
\end{equation}
Each $\check\Phi_j$ is unital and decomposable:
\[
 \check\Phi_j(Y)
  =\nlsum_k R_{j,k}YR_{j,k}^*
   +\nlsum_k T_{j,k}Y^{\mathsf T}T_{j,k}^*,
\]
where the coefficient matrices are real algebraic. On real symmetric
arguments their restrictions admit unital completely positive extensions, and $n_{d,m}$ in the
binomial bound may be replaced by $d+md(d-1)/2$.
In particular, $\check N_{3,m}\le2m^2+1$ and $\check N_{4,m}\le\binom{15m+5}{3m+3}=\exp(O(m))$, while the diagonal estimate remains $\log\check N_{n,n}=O(n^4\log n)$.
\end{prop}
\end{dashedbox}
The bordering construction preserves the factorization
$\Phi_c(aa^*)=U_c(a)U_c(a)^*$ for both types of summands, so the induction
goes through; we count its Kraus operators and kernel degrees in
Appendix~\ref{app:cubic-continuation}, which gives the stated bound.
For $m=1$, the determinant pencil itself has order $d$.

In the table we compare base-ten logarithms of the two bounds for
complex Hermitian inputs: the completely positive column uses
Propositions~\ref{prop:compressed-size} and~\ref{prop:order-three-size},
and the decomposable column uses~\eqref{eq:cubic-seeded-size}.
\begin{center}
\begin{tabular}{rrr}
\toprule
$d=m=n$ & CP bound & Decomposable bound\\
\midrule
3 & 2.117 & 1.279\\
4 & 12.886 & 6.215\\
5 & 109.295 & 43.649\\
6 & 536.237 & 339.120\\
8 & 2528.894 & 1982.235\\
10 & 7770.545 & 6589.210\\
\bottomrule
\end{tabular}
\end{center}
For $d=m=3$ the new pencil has order $19$, and for $d=m=4$ the recurrence
before compression gives order $1\,639\,792$. We also checked the cubic
resolvent identities, the rank-one degree, the signed Kraus
decompositions, and unitality symbolically for $m=2,3,4$.

\section{Pencils for fixed matrices}\label{sec:fixed}
Corollary~\ref{cor:target-intro} concerns the scalar polynomial $\eta(z_1A,\ldots,z_nA)$; we prove it as a special case of Proposition~\ref{prop:psd-tuple}, which allows an arbitrary tuple of positive semidefinite (PSD) matrices scaled by scalar variables. The only point we have to check is that the specialized tuple lies in $\mathcal K_{d,m}^{\circ}$ even when the matrices are singular, and we do so by exhibiting a path to $e_{d,m}$.

\begin{prop}\label{prop:psd-tuple}
Let $B_1,\ldots,B_m\in\HH_d(\F)$ be PSD and suppose
$\eta(B_1,\ldots,B_m)>0$. Put
\[
 f_B(z)=\eta(z_1B_1,\ldots,z_mB_m),\qquad C_j=\Phi_j(B_j),
 \qquad M(z)=\nlsum_jz_jC_j.
\]
Then $f_B$ is hyperbolic with respect to $\one$, $C_j\succeq0$, $M(\one)\succ0$, and, for a polynomial $q$,
\begin{equation*}
 \cone(f_B,\one)=\{z:M(z)\succeq0\},\qquad
 \det M=f_Bq,\qquad \cone(f_B,\one)\subseteq\cone(q,\one).
\end{equation*}
\end{prop}
\begin{proof}
Each $C_j$ is PSD by~\eqref{eq:cp}. To see that the tuple $B=(B_1,\ldots,B_m)$ belongs to $\mathcal K_{d,m}^{\circ}$, we join it to $e_{d,m}$ through $((1-t)B_1+tI_d,(1-t)B_2,\ldots,(1-t)B_m)$ for $0\le t\le1$: all principal minors in every summand are nonnegative, for $t>0$ the summand assigning every index to the first matrix is strictly positive, and the initial value is positive by assumption. Thus $M(\one)=L_{d,m}(B)\succ0$.

Specializing the determinant identity~\eqref{eq:multiplier-main}, we get $\det M=f_Bq$ with $f_B(\one),q(\one)>0$; a definite Hermitian determinant is hyperbolic, so each factor is hyperbolic in that direction. A path in the positive component of $f_B$ starting at $\one$ maps under $z\mapsto(z_1B_1,\ldots,z_mB_m)$ into $\mathcal K_{d,m}^{\circ}$, so $M$ is positive definite along it. Conversely, the set $\{M\succ0\}$ is convex, contains $\one$, and has $f_B>0$ by Theorem~\ref{thm:main}. The open components therefore coincide, and since $M(\one)\succ0$, adding $\eps\one$ gives equality of their closures. Positivity of $q=\det M/f_B$ on the open cone gives the multiplier inclusion.
\end{proof}

\begin{proof}[Proof of Corollary~\ref{cor:target-intro}]
Every coefficient of $f_A$ is nonnegative, and the coefficient obtained from singleton parts is $[z_1\cdots z_n]f_A=n!\nlprod_{i=1}^n A_{ii}$. Hence a positive diagonal gives $f_A(\one)>0$, and we apply Proposition~\ref{prop:psd-tuple} with $B_j=A$. If $A_{ii}=0$, nonnegativity of the $2\times2$ principal minors gives $|A_{ij}|^2\le A_{ii}A_{jj}=0$ for all $j$, so row and column $i$ vanish; each ordered partition then puts $i$ in some part whose principal determinant is zero, and $f_A\equiv0$.

Finally, if $H=\nlsum_jC_j\succ0$, we replace $C_j$ by $H^{-1/2}C_jH^{-1/2}$; the congruence preserves the cone and makes the sum the identity.
\end{proof}

\begin{rmk}[Realification]\label{rmk:realification}
For a Hermitian pencil $M$, whose real part is symmetric and whose imaginary part is skew-symmetric,
\[
 \mathcal R(M)=\begin{pmatrix}\operatorname{Re}M&-\operatorname{Im}M\\\operatorname{Im}M&\operatorname{Re}M\end{pmatrix}
\]
is real symmetric, is PSD exactly when $M$ is PSD, and is definite at the same real directions. Its determinant is $(\det M)^2$. If $\det M=fq$ and the cone of $q$ contains that of $f$, then the realified determinant is $f$ times the multiplier $fq^2$, whose cone still contains that of $f$.

If $M$ has order $N$, we can reduce the real order to $2N-1$.
Every eigenvalue of $\mathcal R(M)$ has even multiplicity, so its negative
eigenspace, if nonzero, meets every real hyperplane, and so does its
kernel whenever $M$ is singular; hence compressing to any fixed real
hyperplane preserves both the positive semidefinite and the positive
definite loci. If $Q$ is an orthonormal basis of that hyperplane and $u$
its unit normal, then
\[
 \det(Q^{\mathsf T}\mathcal R(M)Q)
 =u^{\mathsf T}\operatorname{adj}(\mathcal R(M))u
 =(\det M)u^{\mathsf T}\mathcal R(\operatorname{adj}M)u,
\]
where the second identity follows for invertible $M$ from the inverse
formula and then polynomially. The smaller real pencil therefore retains
the determinant divisibility and the multiplier containment.
\end{rmk}

We give two examples to make the specialization concrete: for rank-one and diagonal matrices the polynomial factors, and at order two the multiplier already appears.

\begin{example}[Rank-one and diagonal matrices]\label{ex:rank-one}
For $A=vv^*$, every principal minor of order at least two vanishes. Only assignments to distinct singleton parts survive, proving~\eqref{eq:rankone-intro}. For diagonal $A$, direct expansion gives $\eta(z_1A,\ldots,z_mA)=\bigl(\nlprod_{i=1}^dA_{ii}\bigr)(z_1+\cdots+z_m)^d$. Singularity of $A$ alone is no obstruction: for $d=m=n$ and $A=\one\one^*$ we get $n!z_1\cdots z_n$, whose cone is the positive orthant.
\end{example}

\begin{example}[Order two]\label{ex:order-two}
For arbitrary Hermitian $2\times2$ inputs, put
\[
 s_i=\nlsum_j(A_j)_{ii},\qquad
 v=((A_1)_{12},\ldots,(A_m)_{12})^{\mathsf T}.
\]
The first induction step gives
\begin{equation*}
 L_{2,m}=\begin{pmatrix}s_1I_m&v\\v^*&s_2\end{pmatrix},\qquad
 \det L_{2,m}=s_1^{m-1}\bigl(s_1s_2-\norm{v}^2\bigr)
              =s_1^{m-1}p_{2,m}.
\end{equation*}
For a fixed PSD matrix $A$, this specializes to
\[
 \eta(z_1A,\ldots,z_mA)
   =A_{11}A_{22}\left(\nlsum_jz_j\right)^2
      -|A_{12}|^2\nlsum_jz_j^2.
\]
The multiplier $s_1^{m-1}$ already occurs at this first nontrivial order.
\end{example}

\subsection{Smaller pencils in the scalar variables}\label{subsec:cubic-scalar-size}
Once the matrices are fixed, only the scalar variables $z$ remain, so we
count the coefficient span of a boundary kernel in $z$ alone, and we can
reduce it twice more: every coordinate of the kernel is divisible by a
fixed monomial in $z$, which lowers its degree, and reducing the
quotient modulo the specialized polynomial removes the monomials
of high degree in one coordinate while keeping a boundary kernel. We
write $\check v_d$ for the polynomial boundary kernel, before compression, of
the pencil obtained in Appendix~\ref{app:cubic-continuation} by continuing
the induction from the cubic pencil of
Proposition~\ref{prop:cubic-resolvent}, and define
\begin{equation}\label{eq:cubic-scalar-degree}
 \check\alpha_{3,m}=0,\qquad
 \check\alpha_{d+1,m}
   =2\check\kappa_{d,m}-1+(m-1)\check\alpha_{d,m},\qquad
 \check E_{d,m}=\check\delta_{d,m}-m\check\alpha_{d,m}.
\end{equation}
Every coordinate of $\check v_d$ has degree at least
$\check\alpha_{d,m}$ in the entries of each argument
(Appendix~\ref{app:cubic-scalar-proofs}), so after the substitution
$A_j\mapsto z_jB_j$ we may divide by $\nlprod_jz_j^{\check\alpha_{d,m}}$. We put
$f_B(z)=\eta(z_1B_1,\ldots,z_mB_m)$. Reducing the quotient
modulo $f_B$ in the coordinate along $\one$ leaves at most
\begin{equation}\label{eq:scalar-monomial-count}
 H_{m,d}(E)=\binom{m+E-1}{m-1}-\binom{m+E-d-1}{m-1}
\end{equation}
monomials at degree $E$, with the second term zero if $E<d$.
For Hermitian inputs we stack the real and imaginary parts of the kernel
before taking its coefficient span, which assigns one real vector to each monomial, so the count is not
doubled.

\begin{dashedbox}[supporting]
\begin{prop}[Real pencils after scalar specialization]\label{prop:cubic-scalar-size}
Let $d\ge3$, $m\ge2$, and let $B_1,\ldots,B_m\in\HH_d(\F)$ be PSD
with $\eta(B_1,\ldots,B_m)>0$. Set
$f_B(z)=\eta(z_1B_1,\ldots,z_mB_m)$, and let $H_{m,d}$ be the monomial count~\eqref{eq:scalar-monomial-count}. There are real symmetric PSD
matrices $C_1,\ldots,C_m$ of order at most $H_{m,d}(\check E_{d,m})$
such that $M_B(z)=\nlsum_jz_jC_j$ is positive definite at $\one$, represents the
open and closed cones of $f_B$, and has determinant $f_Bq$ with
$\cone(f_B,\one)\subseteq\cone(q,\one)$. This real symmetric order
bound holds also for complex Hermitian $B_j$.
\end{prop}
\end{dashedbox}

\begin{dashedbox}[supporting]
\begin{prop}[Cubic bounds uniform in the common matrix]\label{prop:cubic-common-size}
Let $d\ge3$, $m\ge2$, and put
\[
 \check C_{d,m}
  =\binom{m+\check E_{d,m}-1}{m-1}
     \binom{d^2-1+\check\delta_{d,m}}{d^2-1}.
\]
There are unital decomposable maps $\Xi_1,\ldots,\Xi_m$ with fixed
real algebraic coefficient matrices in~\eqref{eq:decomposable-map} such that,
for every $A\in\HH_d(\C)$ with $A\succeq0$ and
$\eta(A,\ldots,A)>0$, the pencil $\nlsum_jz_j\Xi_j(A)$ represents the
open and closed cones of $f_A(z)=\eta(z_1A,\ldots,z_mA)$ and has
determinant $f_Aq$ for a polynomial $q$ with $\cone(f_A,\one)\subseteq\cone(q,\one)$.
The common order of these pencils is at most
\begin{equation}\label{eq:cubic-common-size}
 \min\{\check N_{d,m},\check C_{d,m}\}.
\end{equation}
For real symmetric $A$ alone the maps admit unital completely positive
extensions, and $d^2-1$ in $\check C_{d,m}$ may be replaced by
$d(d+1)/2-1$.

There are also fixed positive unital real-linear maps from
$\HH_d(\C)$ to real symmetric matrices whose pencils have the same
properties for every such Hermitian $A$, with
common order at most
\begin{equation}\label{eq:cubic-common-real-size}
 \min\{2\check N_{d,m}-1,\check C_{d,m}\}.
\end{equation}
Their coefficients in the independent real coordinates of $A$ are
real algebraic.
\end{prop}
\end{dashedbox}
For the uniform bound we count coefficients jointly in $z$ and the
entries of $A$; the perturbation argument that makes the same compression
work for every $A$ is in Appendix~\ref{app:cubic-scalar-proofs}.

At $d=m=n$, these bounds are $\exp(O(n^3\log n))$ for one family valid
for every common matrix and $\exp(O(n^2\log n))$ for a fixed matrix; at
$n=4$ the fixed-matrix bound is $H_{4,4}(11)=244$, which replaces $7\,940$
from the kernel of the completely positive construction (supplement,
Section~S1). For $d=m=3$ the coefficient
count gives six, but the Helton--Vinnikov theorem gives order three for
every such ternary cubic, including singular or reducible
ones~\citep[Theorem~2.2 and Section~4]{helton2007linear}, with PSD
coefficients because the positive orthant lies in the cone.

In the table we compare base-ten logarithms of the bounds; we compared
the universal bounds in Section~\ref{subsec:cubic-seeded-size}. The
common-matrix columns give Hermitian orders, with decomposable maps from
Proposition~\ref{prop:cubic-common-size} and CP maps from the kernel of the
completely positive construction (supplement, Section~S1). For a fixed
matrix, the cubic column gives real symmetric orders even for Hermitian
inputs, and the CP column again comes from the supplement.
\begin{center}
\begin{tabular}{rrrrr}
\toprule
&\multicolumn{2}{c}{Every common $A$}&\multicolumn{2}{c}{One fixed $A$}\\
$n$&CP&Decomposable&CP kernel&Cubic kernel\\
\midrule
4&12.886&6.215&3.900&2.387\\
5&67.765&40.516&9.330&6.440\\
6&147.138&104.086&17.710&13.436\\
8&453.895&369.617&44.070&36.854\\
10&1034.633&894.124&84.282&73.876\\
\bottomrule
\end{tabular}
\end{center}
For rank-one $A$ and $d=m=n$, the cone is the positive orthant and has a
pencil of order $n$ (Example~\ref{ex:rank-one}).
We show in the supplement that the bound $244$ is attained within the
compression method: for an explicit positive definite $4\times4$ matrix,
every constant compression of the uncompressed pencil that represents the
same cone has order at least $244$. For the same matrix, however, a
construction on a common cubic denominator gives a pencil of order at
most $20$, and every real symmetric pencil for that cone has order at
least $5$ (supplement, Section~S2).

\section{Rank-deficient localizations are spectrahedral}\label{sec:localization}
Rank-deficient matrices enter the mixed determinant through limits: as $t\to\infty$, the polynomial $\eta(tA_0,G)$ is dominated by its leading coefficient in $t$, which is again hyperbolic in $G$. Spectrahedrality does not pass to such limits automatically, so we build the representation directly, by compressing the pencil of Theorem~\ref{thm:main} to the kernel of its coefficient at $A_0$. One of these leading coefficients reappears in Appendix~\ref{app:cubic}, where the explicit cubic is a linear section of it.

We fix $A_0\succeq0$ of order $n$ and rank $k<n$ and define
\begin{equation*}
 Q_{A_0}(G)=\sum_{|S|=k}\det A_0[S]\det G[S^c].
\end{equation*}
This polynomial has degree $n-k$, satisfies $Q_{A_0}(I_n)>0$, and is the leading coefficient in
\begin{equation*}
 \eta(tA_0,G)=t^kQ_{A_0}(G)+O(t^{k-1}),\qquad t\longrightarrow+\infty,
\end{equation*}
where the expansion is coefficientwise in $G$. We use the pencil $L_{n,2}(A,B)=\Phi_1(A)+\Phi_2(B)$ of Theorem~\ref{thm:main}, set $H=\Phi_1(A_0)\succeq0$, and let $P$ have orthonormal columns spanning $\ker H$.

\begin{prop}[Compression formula]\label{prop:compression}
The kernel of $H$ is nonzero, the map $\Phi_0(G)=P^*\Phi_2(G)P$ is unital and completely positive in Kraus form, and
\[
 \cone(Q_{A_0},I_n)=\{G:\Phi_0(G)\succeq0\}.
\]
If $\Phi_2(G)=\nlsum_iK_iGK_i^*$, the Kraus operators of $\Phi_0$ are $J_i=P^*K_i$.
\end{prop}
\begin{proof}
Since $\eta(A_0,0)=0$, the open-cone equality of Theorem~\ref{thm:main} forbids $H$ to be positive definite, so $\ker H\ne0$; also $\Phi_0(I_n)=I$.

For any fixed Hermitian $M$ we claim that
\begin{equation}\label{eq:eigen-limit}
 \lambda_{\min}(tH+M)\longrightarrow\lambda_{\min}(P^*MP)
 \quad\text{as }t\longrightarrow+\infty.
\end{equation}
Testing vectors in $\ker H$ gives the upper bound. For the lower bound we choose unit minimizing vectors along a convergent subsequence; the term $t\ip{v}{Hv}$ remains bounded, so the limiting vector lies in $\ker H$, and its $M$-quadratic form gives the reverse inequality in the limit.

By Theorem~\ref{thm:main} and~\eqref{eq:scalar-identities}, which puts $(0,I_n)$ in the same positive component as $e_{n,2}$, the largest zero of $u\mapsto\eta(tA_0,G+uI_n)$ is $-\lambda_{\min}(tH+\Phi_2(G))$. We divide by $t^k$ and let $t\to+\infty$: these real-rooted polynomials converge coefficientwise to $Q_{A_0}(G+uI_n)$, which is nonzero and has leading coefficient $Q_{A_0}(I_n)>0$. Finite zeros of a nonzero polynomial limit are limits of zeros, with multiplicity, for example by the argument principle on small circles around its zeros, so the limit polynomial is real-rooted. By~\eqref{eq:eigen-limit}, the largest zeros of the approximating polynomials have a finite limit; that limit is a zero of the limiting polynomial, and persistence of its other zeros shows that it is the largest zero. The limit equals $-\lambda_{\min}(\Phi_0(G))$, and the root description in Lemma~\ref{lem:component} now gives the asserted cone equality. The Kraus form and unitality follow by compression of $\Phi_2$.
\end{proof}

Adding a zero row and column to $A_0$ and an arbitrary row and column to $G$ gives a polynomial of the same kind. To write it down we use the matrix gradient $C(G)$ of $Q=Q_{A_0}$, with the convention
\begin{equation}\label{eq:matrix-gradient}
 DQ(G)[Y]=\tr(C(G)Y).
\end{equation}
For real symmetric coordinates, this convention means $C_{ii}=\partial Q/\partial G_{ii}$ and $C_{ij}=\tfrac12\partial Q/\partial G_{ij}$ for $i<j$.

\begin{corr}[A zero-border extension]\label{cor:zero-border}
Let
\begin{equation}\label{eq:zero-border}
 Q^\sharp(G,h,s)=sQ(G)-h^*C(G)h.
\end{equation}
With $\Phi_0(G)=\nlsum_iJ_iGJ_i^*$ from Proposition~\ref{prop:compression}, put $R(h)=[J_1h\ \cdots\ J_\ell h]$, where $\ell$ is the number of retained Kraus operators. Then
\begin{equation}\label{eq:zero-border-pencil}
 \cone(Q^\sharp,(I_n,0,1))
 =\left\{(G,h,s):
 \begin{pmatrix}\Phi_0(G)&R(h)\\R(h)^*&sI_\ell\end{pmatrix}\succeq0\right\}.
\end{equation}
In particular, for $G\in\ocone(Q,I_n)$,
\begin{equation}\label{eq:resolvent-identity}
 \frac{h^*C(G)h}{Q(G)}
   =\norm{\Phi_0(G)^{-1/2}R(h)}_{\op}^2.
\end{equation}
\end{corr}
\begin{proof}
We adjoin a zero row and column to $A_0$, obtaining $\widehat A_0=\diag(A_0,0)$. Expanding the principal minors of
$\widehat G=\left(\begin{smallmatrix}G&h\\h^*&s\end{smallmatrix}\right)$ gives
$Q_{\widehat A_0}(\widehat G)=Q^\sharp(G,h,s)$, so Proposition~\ref{prop:compression} gives hyperbolicity with respect to $(I_n,0,1)$.

The polynomial $Q$ is affine along $G\mapsto G-thh^*$, since every principal determinant is, and applying Lemma~\ref{lem:exit} to its compressed pencil we obtain~\eqref{eq:resolvent-identity}. The common-index deletion of $Q^\sharp$ is $Q$, and between the two we need the interlacing that the proof of Lemma~\ref{lem:border-cone} takes from Lemma~\ref{lem:interlacing}. To justify it at the localization, we apply, for $t>0$, the stability of a polynomial plus $v$ times its common-index deletion~\citep[Theorem~2.6]{branden08}, as in the proof of Lemma~\ref{lem:interlacing}, to $t\widehat A_0$ and $\widehat G+uI_{n+1}$, divide by $t^k$, and let $t\to+\infty$. The limit is the nonzero stable polynomial $Q_{\widehat A_0}(\widehat G+uI_{n+1})+vQ(G+uI_n)$, by closure of stability under nonzero coefficient limits, a consequence of Hurwitz's theorem, and Hermite--Biehler gives the interlacing. The component argument of Lemma~\ref{lem:border-cone}, now with one border vector, identifies the open cone by $G\in\ocone(Q,I_n)$ and $s>h^*C(G)h/Q(G)$; by~\eqref{eq:resolvent-identity} and a Schur complement, these two conditions describe exactly the positive definite set in~\eqref{eq:zero-border-pencil}. The pencil equals the identity at $(I_n,0,1)$, so closure gives the full statement.
\end{proof}

\begin{aiusage}
The results in this paper were developed through extensive use of AI assistants. None of the claims were obtained zero-shot: they emerged through sustained back-and-forth discussion, including failed attempts, counterexamples, audits, and repairs. Earlier attempts with OpenAI's GPT~Astra (Pro) to extend Br\"and\'en's matrix--tree argument~\citep{branden2014hyperbolicity} to the fixed-matrix polynomial $\eta(z_1A,\ldots,z_nA)$ did not produce a proof. The investigation then searched for a counterexample, adapting the norm obstruction of~\citet{openai2026nonspectrahedral} to the cubic of Appendix~\ref{app:cubic}; restricted Clifford and exceptional-Jordan-algebra ansatzes and a sextic square-domination candidate were also explored. The obstruction proved valid for pure powers but not for cones, and the author's suggestion to regard the mixed determinant as an ``overgrown elementary symmetric polynomial'' turned the search back to the positive direction, where GPT~Astra (Pro) found the construction of Theorem~\ref{thm:main}. The search contributed the decision to permit determinant multipliers, to use extremal eigenvalues rather than traces, and to preserve linear rank-one factorizations at every step of the induction. 

The counterexamples of Appendix~\ref{sec:aim} come from an earlier, separate investigation~\citep{Sra2026CountExMachina}. Those for Problems~35 and~38 were found with GPT-5~Pro in January 2026 and those for Problems~36 and~37 with GPT-5.6~Pro in July 2026; Anthropic's Claude (Opus~5) isolated the Schur-positive witness for Problem~37 in August 2026, and the author formalized and audited all four. For this version, Claude (Opus~5.5) recomputed every expansion, minor, and inequality of Appendix~\ref{sec:aim} in exact arithmetic. The finite exact checks mentioned in Section~\ref{subsec:sizes} and Appendix~\ref{app:cubic} were run within the same AI-assisted workflow and were not rerun for this version. Anthropic's Claude (Opus~5.5) prepared the manuscript before the compression bounds were added, and Anthropic's Claude (Fable~5.1 and~5.2) revised the presentation of the whole manuscript in this version. 

OpenAI's Codex developed and integrated the bounds in Section~\ref{sec:smaller}, and the specialization bounds in Appendix~\ref{app:scalar-proofs} and in the supplement, including the cubic resolvent construction and its decomposable continuation, independent assistant audits, and exact checks of the polynomial-kernel identities, degree bounds, coefficient ranks, and signed Kraus decompositions. A subsequent dimension-optimization package supplied the additional scalar bounds and the sharp constant-compression example in Section~\ref{subsec:cubic-scalar-size} and in the supplement. Its exact checks of the cubic pencil and of the induction, and its finite-field rank, smoothness, and corank certificates were rerun during integration. Assistant checks are not external peer review or formal verification.
\end{aiusage}

\aiacknowledgements

\appendix
\addtocontents{toc}{\protect\setcounter{tocdepth}{1}}
\begingroup
\renewcommand{\contentsname}{Appendix contents}
\tableofcontents
\endgroup
\section{Proofs of the size bounds for the universal cone}\label{app:kernel-counts}
The compressed bounds of Section~\ref{sec:smaller} depend on the degrees
of polynomial boundary kernels. In the constructions below
we keep track of the rank-one degree, which we need to clear denominators,
and of the nonvanishing, which we need for the cone equality.

\subsection{The kernel for the completely positive construction}\label{app:kernel-construction}
\begin{proof}[Proof of Proposition~\ref{prop:compressed-size}]
We construct homogeneous polynomial vectors $v_d,g_d$ such that
\begin{equation}\label{eq:kernel-polynomial}
 L_{d,m}v_d=p_{d,m}g_d,
\end{equation}
with $\deg v_d=\delta_{d,m}$ and $v_d$ nonzero on a boundary chart. We also require that, for every rank-one direction $h$ in any one argument,
\begin{equation}\label{eq:rank-one-kernel-degree}
 \deg_t v_d(X+th)\le\kappa_{d,m},\qquad
 \deg_t g_d(X+th)\le\kappa_{d,m}-1.
\end{equation}
We work in the coordinates given by the diagonal sums $s_i=\nlsum_j(A_j)_{ii}$ and the ordered off-diagonal entries; by the partition formula and the block recursion, the polynomial and the pencil depend on the diagonal entries only through these sums. In the Hermitian case an entry and its conjugate count as independent formal variables; the vectors we construct have integer coefficients.

At order two we put $b=((A_1)_{12},\ldots,(A_m)_{12})^{\mathsf T}$. Then
\[
 L_{2,m}=\begin{pmatrix}s_1I_m&b\\b^*&s_2\end{pmatrix},\qquad
 v_2=\begin{pmatrix}-b\\s_1\end{pmatrix},\qquad
 L_{2,m}v_2=\begin{pmatrix}0\\p_{2,m}\end{pmatrix}.
\]
This vector has degree one, is independent of $s_2$, and is nonzero on the boundary chart $s_1>0$; the vector $g_2$ is constant, so~\eqref{eq:rank-one-kernel-degree} holds with $\kappa_{2,m}=1$.

Suppose that these properties hold at order $d$, and write $\delta=\delta_{d,m}$ and $\kappa=\kappa_{d,m}$. For an interior base tuple $X$ and border vectors $a_j$ we put
\[
 p=p_{d,m}(X),\qquad D_j=D_{a_ja_j^*}^{(j)}p(X),\qquad
 Y_j=X-\frac{p}{D_j}h_j,
\]
where $h_j$ has $a_ja_j^*$ in argument $j$ and zeros elsewhere; we assume $D_j\ne0$ for the moment. By rank-one affinity, $p_{d,m}(Y_j)=0$. We define
\[
 w_j=D_j^\kappa v_d(Y_j),\qquad
 z_j=\frac{U_j(a_j)^*w_j}{D_j}.
\]
The first degree bound in~\eqref{eq:rank-one-kernel-degree} makes $w_j$ polynomial, and we claim that the division defining $z_j$ is exact. Indeed, write $v_d(X+th_j)=\nlsum_{\ell=0}^{\kappa}t^\ell V_\ell(X,h_j)$. Since $p_{d,m}$ is affine along $h_j$, the coefficient of $t^{\kappa+1}$ in~\eqref{eq:kernel-polynomial} is
$U_jU_j^*V_\kappa=0$; positivity gives $U_j^*V_\kappa=0$ on Hermitian inputs, hence as a polynomial identity, and in
$U_j^*w_j=\nlsum_\ell D_j^{\kappa-\ell}(-p)^\ell U_j^*V_\ell$ the only term without a factor $D_j$ vanishes. The degrees of $w_j$ and $z_j$ are
$\delta+(d+1)\kappa$ and $\delta+(d+1)\kappa-d$, respectively. Substituting in~\eqref{eq:kernel-polynomial} gives
\begin{equation}\label{eq:polynomial-kernel-identity}
 L(X)w_j=pU_j(a_j)z_j,\qquad U_j(a_j)^*w_j=D_jz_j.
\end{equation}
Identity~\eqref{eq:polynomial-kernel-identity} is polynomial and therefore also holds when $D_j=0$.

On the spaces of~\eqref{eq:new-pencil} we set $Z=z_1\otimes\cdots\otimes z_m$ and
\begin{equation}\label{eq:next-kernel}
 (v_{d+1})_T=pZ,\qquad
 (v_{d+1})_{V_j}
 =-z_1\otimes\cdots\otimes w_j\otimes\cdots\otimes z_m.
\end{equation}
By~\eqref{eq:polynomial-kernel-identity}, every top block of $L^+v_{d+1}$ cancels, and the bottom block is $(\sigma p-\nlsum_jD_j)Z=p_{d+1,m}Z$. Thus~\eqref{eq:kernel-polynomial} holds at order $d+1$, with $g_{d+1}=(0,\ldots,0,Z)$. Every block in~\eqref{eq:next-kernel} has degree
\begin{equation}\label{eq:kernel-degree-step}
 \delta_{d+1,m}=m\bigl[\delta+(d+1)\kappa\bigr]-d(m-1),
\end{equation}
and the vector is independent of the new diagonal sum $s_{d+1}=\sigma$.

To verify~\eqref{eq:rank-one-kernel-degree} at the next order, we update the bordered argument $A_c^+$ in direction $(b,\beta)(b,\beta)^*$: the base becomes $X(t)=X+th_c(b)$, $a_c(t)=a_c+tb\overline\beta$, and the other borders remain fixed. For $j\ne c$, the polynomial $D_j(t)$ is generically linear, and $p(X(t))/D_j(t)$ is bounded at infinity; applying the uniform rank-one degree bound to
\[
 Y_j(t)=X+th_c(b)-\frac{p(X(t))}{D_j(t)}h_j(a_j)
\]
we get $v_d(Y_j(t))=O(t^\kappa)$, hence $\deg_t w_j\le2\kappa$ and $\deg_t z_j\le2\kappa-1$.

For $j=c$ we put $p_1=D_{bb^*}^{(c)}p(X)$. The polynomial $D_c(t)$ is generically quadratic, with leading coefficient $|\beta|^2p_1$; to verify this coefficient we polarize the identity $\bigl(D_{aa^*}^{(c)}\bigr)^2p=0$, which expresses rank-one affinity: the coefficient of $s^3$ in $\bigl(D_{(a+sb)(a+sb)^*}^{(c)}\bigr)^2p=0$ gives $D_{bb^*}^{(c)}D_{ba^*+ab^*}^{(c)}p=0$, and together with $\bigl(D_{bb^*}^{(c)}\bigr)^2p=0$, this identity gives the asserted leading coefficient after we absorb $\beta$ into the cross term. Consequently, the leading $tbb^*$ terms in
\[
 Y_c(t)=X(t)-\frac{p(X(t))}{D_c(t)}h_c(a_c(t))
\]
cancel, so $Y_c(t)$ is bounded at infinity, hence $\deg_t w_c\le2\kappa$; since $U_c(a_c(t))$ is linear and $D_c(t)$ is quadratic, $\deg_t z_c\le2\kappa-1$ as well. The leading coefficients we used are nonzero on open parameter sets: for $j\ne c$ we take two distinct coordinate rank-one directions at $e_{d,m}$, and for $j=c$ we take $\beta p_1\ne0$. As all the vectors are polynomial, these degree bounds hold identically.

The bottom block $pZ$ and every top block in~\eqref{eq:next-kernel} therefore have rank-one degree at most $\kappa_{d+1,m}=m(2\kappa_{d,m}-1)+1$, whereas $g_{d+1}$ has degree at most $\kappa_{d+1,m}-1$. Starting from $\kappa_{2,m}=\delta_{2,m}=1$, this recurrence and~\eqref{eq:kernel-degree-step} give exactly~\eqref{eq:kernel-degree}.

For nonvanishing we fix an argument $j$ and a coordinate rank-one direction $h_j$, whose directional derivative has value one at $e_{d,m}$. In linear coordinates $(t,\xi)$ along $h_j$, rank-one affinity gives $p_{d,m}(\xi+th_j)=a(\xi)t+b(\xi)$, and on $a(\xi)\ne0$ the exit map sends an interior tuple to $\xi-h_jb(\xi)/a(\xi)$. Interior tuples near $e_{d,m}$ project onto an open set of $\xi$, so their exits are Zariski dense in the irreducible hypersurface $\{p_{d,m}=0\}$; since $v_d$ does not vanish identically there, $w_j$ is a nonzero polynomial vector. At interior points,~\eqref{eq:polynomial-kernel-identity} shows that $z_j=0$ implies $w_j=0$. The conditions $D_j\ne0$ and $z_j\ne0$ for all $j$ therefore hold together on an open dense set of the base and border parameters. There we set $\sigma=\nlsum_jD_j/p$; Lemma~\ref{lem:border-cone} and closure place the tuple on the boundary, where the derivative of $p_{d+1,m}$ in $\sigma$ is $p>0$ and $(v_{d+1})_T=pZ\ne0$. The induction is complete.

We have $n_{d,m}$ coordinates, and $v_d$ is independent of $s_d$, so its coefficient span has dimension at most the binomial term in~\eqref{eq:compressed-size}. We choose a real orthonormal basis $Q$ of the span of its integer coefficient vectors; for Hermitian inputs, evaluations of $v_d$ lie in the complex column space of the same real matrix $Q$, so Lemma~\ref{lem:kernel-compression} and Proposition~\ref{prop:irreducible} give a compression over both fields. Rational linear algebra followed by orthonormalization makes $Q$ real algebraic, so the compressed Kraus operators are real algebraic as well. Identifying the two ordered entries of each off-diagonal pair before taking the span gives the smaller count for real symmetric inputs alone.

The determinant quotient has value one at $e_{d,m}$, is hyperbolic, and has a cone containing $\mathcal K_{d,m}$ by the argument of Proposition~\ref{prop:multiplier}. If the original pencil is smaller, we retain it instead, which gives the minimum in~\eqref{eq:compressed-size}. To bound the number of Kraus operators, we form $C_j=\nlsum_k\operatorname{vec}(\widehat K_{j,k})\operatorname{vec}(\widehat K_{j,k})^{\mathsf T}$; factoring this PSD matrix of order $d\widehat N_{d,m}$ and reshaping its factor columns gives at most $d\widehat N_{d,m}$ real algebraic Kraus operators for the same map.

Finally, $\kappa_{d,m}\le(2m)^{d-2}$, so
$\log(\delta_{d,m}+1)=O(d\log m+\log d)$ for $d,m\ge2$, and combining this
estimate with $n_{d,m}=O(md^2)$ and
$\binom{a+b}{a}\le(e(a+b)/a)^a$ we obtain the asserted bound.
\end{proof}

\subsection{An order-three bound with completely positive maps}\label{subsec:order-three-size}
At order three we can group the kernel coordinates and so obtain a bound that retains completely positive maps over both fields, without computing the rank of their span.

\begin{prop}\label{prop:order-three-size}
Let $m\ge2$, and let $\mathrm{Del}_m=\nlsum_{k=0}^m\binom mk\binom{m+k}{k}$ be the central Delannoy number. The cone $\mathcal K_{3,m}$ has a pencil with the properties in Proposition~\ref{prop:compressed-size} and order at most
\begin{equation}\label{eq:delannoy-bound}
 \widetilde N_{3,m}
 =\frac{(m+1)\mathrm{Del}_m+(3m+1)\mathrm{Del}_{m-1}}2
       -m(3m+1)2^{m-2}.
\end{equation}
\end{prop}
\begin{proof}
We write
\[
 L_{2,m}=\begin{pmatrix}sI_m&b\\b^*&t\end{pmatrix},\qquad
 p=st-b^*b,\qquad a_j=(u_j,v_j)^{\mathsf T},
\]
and set
\begin{align*}
 q_j&=t|u_j|^2+s|v_j|^2-u_j\overline v_j\overline b_j-\overline u_jv_jb_j,\\
 \zeta_j&=s\overline v_j e_j-\overline u_jb,\qquad
 c_j=\begin{pmatrix}p u_j\overline v_j e_j-q_jb\\s q_j-p|u_j|^2\end{pmatrix},
\end{align*}
where $e_j$ is the $j$th coordinate vector in $\F^m$. Direct multiplication gives
\[
 L_{2,m}c_j=pU_j(a_j)\zeta_j,\qquad
 U_j(a_j)^*c_j=q_j\zeta_j,
\]
so the vector with bottom block $p\bigotimes_j\zeta_j$ and $j$th top block
$-\zeta_1\otimes\cdots\otimes c_j\otimes\cdots\otimes\zeta_m$ satisfies~\eqref{eq:kernel-polynomial}; its degree is $2m+2$, and at the smooth boundary point $s=t=1$, $b=0$, $u_j=0$, $v_j=1$, and $\sigma=m$ its bottom block is nonzero.

We now count groups of identical coordinates. For a bottom coordinate indexed by a word $(i_1,\ldots,i_m)$, we fix the positions $F=\{j:i_j=j\}$; the coordinate polynomial depends only on $F$ and the multiplicities of the remaining letters. If $k=m-|F|$, all degree-$k$ multiplicity vectors occur except the $k$ monochromatic ones in a forbidden color: we assign the prescribed copies of the letters to the remaining positions, each of which forbids its own color, and since any two positions together allow every letter, Hall's condition reduces to the singleton conditions, which exclude exactly those monochromatic choices. Thus we have at most
\[
 \beta_m=\nlsum_{k=0}^m\binom mk
       \left[\binom{m+k-1}{k}-k\right]
       =B_m-m2^{m-1},\qquad
 B_m=\nlsum_{k=0}^m\binom mk\binom{m+k-1}{k}
\]
groups of identical bottom coordinates.

In top block $j$, the tensor positions $\ell\ne j$ are \emph{spectators}, and we fix those with $i_\ell=\ell$. If $k$ spectator positions remain, the coordinate $i_j=j$ and the last coordinate of $c_j$ each give $\binom{m+k-1}{k}-k$ groups; for the other coordinates, $(c_j)_{i_j}=-q_jb_{i_j}$ adds a letter forbidden to equal $j$, and the same count for these $k+1$ positions gives $\binom{m+k}{k+1}-(k+1)$ groups. Each top block thus has at most
\begin{align*}
 \omega_m&=\nlsum_{k=0}^{m-1}\binom{m-1}{k}
       \left[2\binom{m+k-1}{k}+\binom{m+k}{k+1}-3k-1\right]\\
    &=B_m+\mathrm{Del}_{m-1}-(3m-1)2^{m-2}
\end{align*}
groups. Their indicator vectors span a fixed real space containing all coefficient vectors, so Lemma~\ref{lem:kernel-compression} gives a pencil of order at most $m\omega_m+\beta_m$, and the identity $B_m=(\mathrm{Del}_m+\mathrm{Del}_{m-1})/2$ yields~\eqref{eq:delannoy-bound}.
\end{proof}

Since $\binom{m+k}{k}\le2^{m+k}$, we have $\mathrm{Del}_m\le6^m$, and hence $\widetilde N_{3,m}=O(m6^m)$. The first values are
\begin{center}
\begin{tabular}{rrrrrr}
\toprule
$m$&2&3&4&5&6\\
\midrule
$N_{3,m}=(m+2)m^m$&16&135&1536&21875&373248\\
$\widetilde N_{3,m}$&16&131&1004&6977&45626\\
\bottomrule
\end{tabular}
\end{center}
For $m=3,4$ we computed the exact rational coefficient ranks of the vector constructed in the proof, $131$ and $1004$, both for real symmetric inputs and with conjugate coordinates treated formally as independent, and we checked the polynomial identity~\eqref{eq:kernel-polynomial} and the identities~\eqref{eq:polynomial-kernel-identity} in exact rational arithmetic for $(d,m)=(3,2),(3,3),(3,4),(4,2)$ over both fields.

\subsection{Continuation from the cubic pencil}\label{app:cubic-continuation}
\begin{proof}[Proof of Proposition~\ref{prop:cubic-seeded-size}]
We first check that the bordering induction preserves the stated class
of maps. We write every old summand as
$K_{c,k}\tau_{\varepsilon_{c,k}}(Y)K_{c,k}^*$, where
$\tau_0(Y)=Y$ and $\tau_1(Y)=Y^{\mathsf T}$, and pad with zero summands
so that each argument has $r$ Kraus operators. On a rank-one input $aa^*$
we form $U_c(a)$ with column $K_{c,k}a$ when $\varepsilon_{c,k}=0$ and
$K_{c,k}\overline a$ when $\varepsilon_{c,k}=1$; then $\Phi_c(aa^*)=U_c(a)U_c(a)^*$, so the exit-threshold and tensor-sum arguments apply unchanged.

We use the embeddings $K^+_{c,j,\mathbf k}$ of
\eqref{eq:newK-top}--\eqref{eq:newK-bottom}, assigning each the type
$\varepsilon_{c,k_j}$, and for type one we apply the embedding to the
transpose of the bordered argument $A_c^+$, whose border is $\overline{a_c}$
while its new diagonal entry is still $\alpha_c$. Thus the
top-to-bottom block is exactly the column of $U_c(a_c)$ required by
\eqref{eq:new-pencil}, the top diagonal blocks remain copies of the
old maps, and the bottom diagonal block is $\alpha_cI$; the calculation
proving~\eqref{eq:unital-closes} is unchanged. The class is therefore
closed under bordering, with real coefficient matrices at every stage,
and over real symmetric inputs we may treat all the transposed summands
as ordinary completely positive summands.

Starting from Proposition~\ref{prop:cubic-resolvent}, the order and
the number of Kraus operators therefore satisfy
\[
 S_{3,m}=2m^2+1,\qquad \rho_{3,m}=2m^2-3m+2,
\]
\[
 S_{d+1,m}=mS_{d,m}\rho_{d,m}^{m-1}+\rho_{d,m}^m,
 \qquad \rho_{d+1,m}=m\rho_{d,m}^m.
\]
We solve the second recurrence and then use
$S_{d+1,m}/\rho_{d+1,m}=S_{d,m}/\rho_{d,m}+1/m$ to get the stated
formulas, which describe the pencils before coefficient compression.

The kernel~\eqref{eq:cubic-resolvent-kernel} of the cubic pencil has
total degree two, rank-one degree one, and constant right-hand side
$e_{\mathrm{last}}$. In the proof of
Proposition~\ref{prop:compressed-size} we used only the factorization
$\Phi_c(aa^*)=U_c(a)U_c(a)^*$, positivity, and polynomial dependence
in the real coordinates, and all three hold here; in particular, the
highest rank-one coefficient still satisfies
$U_cU_c^*V_\kappa=0$, hence $U_c^*V_\kappa=0$, so the division by
$D_c$ remains exact, and the vectors $U_c(a_c(t))$ are affine in the
real update parameter even for transposed summands. The same degree
calculation therefore gives
\[
 \check\kappa_{3,m}=1,\qquad
 \check\kappa_{d+1,m}=2m\check\kappa_{d,m}-m+1,
 \qquad
 \check\delta_{d,m}=(d-1)\check\kappa_{d,m}.
\]
The generic nonvanishing argument still holds: rank-one exits
are dense in the irreducible boundary hypersurface, and at interior
base points the identity $Lw_c=pU_cz_c$ shows that $z_c=0$ forces
$w_c=0$, so the tensor product in the next boundary kernel
is nonzero on the required open chart.

The cubic pencil and its kernel depend on diagonal entries only through
their sums; bordering preserves this property, and each new kernel
is independent of the last diagonal sum. In the ordered-entry
coordinates, transposition simply exchanges two entries and introduces
no nonreal scalar coefficients, so the kernel has real algebraic
coefficient vectors in $n_{d,m}-1$ variables, whose span has dimension
at most the binomial term in~\eqref{eq:cubic-seeded-size}.
Lemma~\ref{lem:kernel-compression} then gives the cone equalities and
determinant divisibility; we compress the two types of summands
separately, which preserves unital decomposability, and a real
algebraic orthonormal basis works over both fields. The determinant
quotient is hyperbolic and its cone contains the original cone, by the
argument of Proposition~\ref{prop:multiplier}.

Retaining the uncompressed pencil or the earlier completely positive
pencil gives the other two terms of the minimum, and identifying
opposite off-diagonal entries gives the stated real symmetric variable
count. Finally, $\check\delta_{4,m}=3m+3$, which gives the displayed
order-four bound, and the diagonal asymptotic follows by the same
binomial estimate as in Proposition~\ref{prop:compressed-size}.
\end{proof}

\section{Proofs of the size bounds for fixed matrices}\label{app:scalar-proofs}
Restricted to scalar multiples of fixed matrices, a polynomial boundary kernel can vanish identically, and the specialized polynomial may become reducible. We use the kernel continued from the cubic pencil in Appendix~\ref{app:cubic-continuation} and the exit identities of Appendix~\ref{app:kernel-counts}, and we retain the determinant factor through every specialization by generic coprimality and positive definite perturbations.

\subsection{The degree in each argument}\label{app:cubic-scalar-proofs}
We write $\check v_d$ for the kernel of Appendix~\ref{app:cubic-continuation} and $\check\kappa_{d,m}$ for its rank-one degree. Every coordinate of the cubic kernel~\eqref{eq:cubic-resolvent-kernel} has degree at least $\check\alpha_{3,m}=0$ in the entries of each argument, and scaling one argument gives the recurrence~\eqref{eq:cubic-scalar-degree}. To check the inductive step, we scale the bordered argument $A_c^+$ by $t$ and work at $t=0$. On a nonempty generic chart, $D_c$ has valuation two and $Y_c$ remains bounded, so $w_c$ and $z_c=U_c^*w_c/D_c$ have valuations at least $2\check\kappa_{d,m}$ and $2\check\kappa_{d,m}-1$; transposed summands still satisfy $U_c(ta)=tU_c(a)$ for real $t$. For $j\ne c$, $D_j(0)\ne0$ and the $c$th argument of $Y_j$ is still $tA_c$, so both $w_j$ and $z_j$ have valuation at least $\check\alpha_{d,m}$. The tensor formulas~\eqref{eq:next-kernel} give the recurrence, and since the degree bounds are polynomial identities they extend off the chart. In particular, $\check v_d(z_1B_1,\ldots,z_mB_m)$ is divisible by $\nlprod_jz_j^{\check\alpha_{d,m}}$.

\begin{proof}[Proof of Proposition~\ref{prop:cubic-scalar-size}]
We write $p=p_{d,m}$, $\delta=\check\delta_{d,m}$, $\alpha=\check\alpha_{d,m}$, and $E=\check E_{d,m}$, and let $L$ and $v$ denote the uncompressed pencil and polynomial boundary kernel continued from the cubic pencil of Proposition~\ref{prop:cubic-resolvent}. We treat the entries of $B_1,\ldots,B_m$ as independent variables and put $c(z)=\nlprod_jz_j^\alpha$. The divisibility just proved gives us polynomial vectors
\[
 \widetilde v(B,z)=c(z)^{-1}v(z_1B_1,\ldots,z_mB_m),\qquad
 L(z_1B_1,\ldots,z_mB_m)\widetilde v=f_B\widetilde g.
\]
Here $\widetilde g$ is polynomial as well: with $B$ still indeterminate, $f_B$ has no $z_j$ factor when $m\ge2$, so the original kernel identity implies that $c$ divides the specialized right-hand vector. The vector $\widetilde v$ has degree $E$ in $z$ and degree $\delta$ in $B$.

We first justify a generic nonvanishing condition within this family of sections. Proposition~\ref{prop:irreducible} and its proof remain valid over $\C$. We choose a coordinate $h$ of $v$ that is nonzero on a real symmetric boundary chart, so $p\nmid h$, and let $\widetilde h=c^{-1}h(z_1B_1,\ldots,z_mB_m)$. After inverting all the $z_j$, the change of variables $X_j=z_jB_j$ identifies $f_B$ with the irreducible polynomial $p(X)$ and shows that $f_B$ does not divide $\widetilde h$. Without inverting the $z_j$, any additional factor of $f_B$ would have to be a monomial in the $z_j$, which is impossible because $p$ remains nonzero when any one argument is zero. Thus $f_B$ is irreducible, and coprime to $\widetilde h$, also over the fraction field $\C(B)$ by Gauss's lemma.

Hence the resultant
\[
 R(B,u)=\operatorname{Res}_t\bigl(f_B(u+t\one),
                       \widetilde h(B,u+t\one)\bigr)
\]
is not identically zero, and the same argument applies to real symmetric indeterminates. We may therefore choose $u$ and real symmetric positive definite matrices $\Delta_1,\ldots,\Delta_m$, all with rational coordinates, such that $R(\Delta,u)\ne0$: these tuples form a nonempty real open set. The resultant is homogeneous in $B$, since its two inputs have respective $B$-degrees $d$ and $\delta$ and generic $t$-degrees $d$ and $E$. Consequently $R(B+\varepsilon\Delta,u)$ is a nonzero polynomial in $\varepsilon$, with nonzero leading coefficient $R(\Delta,u)$, and it is nonzero for all sufficiently small positive $\varepsilon$.

We temporarily write $B$ for one of these perturbed tuples. Then $f_B$ and $\widetilde h(B,\cdot)$ are coprime in $z$: any nonconstant homogeneous common factor $\phi$ satisfies $\phi(\one)\ne0$, because $f_B(\one)>0$, so $\phi(u+t\one)$ has degree $\deg\phi>0$, contradicting the nonzero line resultant; this argument does not require $f_B$ to be irreducible. We divide each coordinate of $\widetilde v$ by $f_B$ as a polynomial in the coordinate along $\one$, whose leading coefficient $f_B(\one)>0$ is constant in the remaining coordinates. The resulting remainder $r_B$ is homogeneous of degree $E$, has degree less than $d$ in that coordinate, has at most $H=H_{m,d}(E)$ monomials, and satisfies $L(z_1B_1,\ldots,z_mB_m)r_B=f_B\gamma_B$ for a polynomial vector $\gamma_B$; its coordinate obtained from $\widetilde h$ remains coprime to $f_B$.

We now use the realification $\mathcal R$ of Remark~\ref{rmk:realification} and put
\[
 \rho(r_B)=\binom{\operatorname{Re}r_B}{\operatorname{Im}r_B}.
\]
For real scalar variables the kernel identity becomes $\mathcal R\bigl(L(z_1B_1,\ldots,z_mB_m)\bigr)\rho(r_B)=f_B\rho(\gamma_B)$, and if $r_B=\nlsum_\beta r_\beta z^\beta$ then $\rho(r_B)=\nlsum_\beta\rho(r_\beta)z^\beta$, so the real coefficient span has dimension at most $H$. We choose a real orthonormal basis $Q$ of this span and compress the realified pencil to obtain $M_B$, whose coefficients are PSD and which satisfies $M_B(\one)\succ0$.

We write the coprime coordinate as $h=a+ib$, where $a$ and $b$ are real coordinate polynomials of $\rho(r_B)$. Multiplying the compressed kernel identity by $\operatorname{adj}M_B$ shows that $f_B$ divides both $(\det M_B)a$ and $(\det M_B)b$, hence $(\det M_B)h$, and coprimality gives $f_B\mid\det M_B$. The positive definite locus of $M_B$ is convex, contains $\one$, and avoids $\{f_B=0\}$; compression gives the reverse inclusion from the open cone of $f_B$, and taking closures gives the closed-cone equality. Positivity of $\det M_B/f_B$ on the open cone gives the multiplier inclusion.

We now return to the original tuple and let $\varepsilon$ tend to zero, passing to a subsequence of the real isometries with constant order and limit $Q$. Strict feasibility persists because the original tuple belongs to $\mathcal K_{d,m}^{\circ}$, so $L(B_1,\ldots,B_m)\succ0$. The determinant factorization persists too: if $\det M_\varepsilon=f_\varepsilon q_\varepsilon$, multiplication by the nonzero limit polynomial $f_B$ is an injective linear map between fixed finite-dimensional coefficient spaces; its smallest singular value is positive and stays bounded below for nearby polynomials, so the coefficients of $q_\varepsilon$ are bounded and a subsequence converges. Thus $f_B\mid\det M_B$, and the same positive-component argument gives both cone equalities at the limit; realification and compression preserve the positivity of the coefficients.
\end{proof}

\begin{proof}[Proof of Proposition~\ref{prop:cubic-common-size}]
We keep the notation of the preceding proof and choose real symmetric positive definite matrices $\Delta_1,\ldots,\Delta_m$ with rational entries for which its generic nonvanishing argument applies, with $h$ the chosen coordinate of $\check v_d$ and $\widetilde h$ its divided restriction. At the tuple $\Delta$, $f_\Delta$ has no coordinate factor $z_j$: setting $z_j=0$ leaves the nonzero term $\det(\Delta_k)z_k^d$ for any $k\ne j$. Hence coprimality with the divided coordinate $\widetilde h$ also gives coprimality with the original restricted coordinate $h(z_1\Delta_1,\ldots,z_m\Delta_m)=c(z)\widetilde h(\Delta,z)$. For $\varepsilon>0$ we put
\[
 r_{A,\varepsilon}(z)
 =\frac{\check v_d\bigl(z_1(A+\varepsilon \Delta_1),\ldots,
                  z_m(A+\varepsilon \Delta_m)\bigr)}
        {\nlprod_jz_j^{\check\alpha_{d,m}}}.
\]
The degree bound in each argument makes $r_{A,\varepsilon}$ polynomial. It has degree $\check E_{d,m}$ in $z$ and degree $\check\delta_{d,m}$ jointly in the entries of $A$ and $\varepsilon$, so its coefficient of $A^\mu z^\beta$ is a fixed vector times $\varepsilon^{\check\delta_{d,m}-|\mu|}$, and the span $W$ of these coefficient vectors is independent of $\varepsilon\ne0$. Since $\check v_d$ omits the last diagonal sum, $A_{dd}$ does not occur. We count the degree-$\check E_{d,m}$ monomials in $z$ and the monomials of degree at most $\check\delta_{d,m}$ in the remaining $d^2-1$ ordered-entry coordinates of $A$, which gives $\check C_{d,m}$; in the Hermitian case we count the ordered off-diagonal entries as formal variables, as in the universal compression argument, and for real symmetric $A$ we identify opposite off-diagonal entries, which leaves $d(d+1)/2-1$ coordinates. The coefficient vectors are real algebraic, and a real algebraic orthonormal basis $Q$ of $W$ preserves the completely positive and the transposed completely positive summands separately. We set $\Xi_j(Y)=Q^*\check\Phi_j(Y)Q$.

For each fixed $\varepsilon>0$, the restrictions of $p_{d,m}$ and $h$ along $z\mapsto(z_1(A+\varepsilon \Delta_1),\ldots,z_m(A+\varepsilon \Delta_m))$ are coprime for generic $A$, because a line resultant witnessing coprimality is nonzero at $A=0$, and dividing the restricted kernel by $\nlprod_jz_j^{\check\alpha_{d,m}}$ preserves coprimality and the polynomial kernel identity. The coefficient compression argument of the preceding proof therefore gives us determinant divisibility and the cone equalities for generic positive definite $A$, and its limiting argument extends them to all $A\succeq0$, since the uncompressed pencil satisfies $L(z_1(A+\varepsilon \Delta_1),\ldots)\succeq\varepsilon\bigl(\nlsum_j\lambda_{\min}(\Delta_j)\bigr)I$ at $z=\one$. Finally we let $\varepsilon\downarrow0$. We use the same $Q$ throughout, and $\eta(A,\ldots,A)>0$ ensures strict feasibility at $\one$ in the limit; divisibility and the cone equalities pass to the limit as before, and positivity and separate unitality of the decomposable maps persist. Retaining the universal pencil gives the other term in~\eqref{eq:cubic-common-size}.

For a real pencil valid for all Hermitian $A$, we instead express the entries of $A$ in independent real coordinates and use $\rho(r_{A,\varepsilon})$ with the realified pencil. There are still $d^2-1$ coordinates, and each monomial contributes one real coefficient vector, so the count is again $\check C_{d,m}$. The joint real and imaginary coprimality argument of the preceding proof and the same perturbation limits give a real compression independent of $A$, and realification and compression preserve positivity and unitality. For the remaining bound in~\eqref{eq:cubic-common-real-size} we apply Remark~\ref{rmk:realification} to the universal pencil: its fixed real hyperplane compression has order $2\check N_{d,m}-1$ and retains the cone and multiplier properties.
\end{proof}

\section{Counterexamples to the Borcea--Br\"and\'en problems}\label{sec:aim}
Corollary~\ref{cor:target-intro} represents the cone of $f_A(z)=\eta(z_1A,\ldots,z_nA)$ for every positive semidefinite $A$ with positive diagonal; in this appendix $A\in\HH_d(\F)$ may have any order $d$, and Proposition~\ref{prop:psd-tuple} with $B_1=\cdots=B_n=A$ represents the cone whenever $f_A(\one)>0$. Problem~35 of \citet{aim2007polyaschurlax} asks whether every homogeneous symmetric real stable polynomial is of the form $f_A$, and three companion problems from the same list concern the Schur expansions of such polynomials. In this appendix we answer all four negatively. The counterexamples were found with AI assistance and are archived, with certificates in exact arithmetic, as the case \texttt{aim-problems} of~\citet{Sra2026CountExMachina}.

The problems concern Schur expansions, and for the Schur comparisons we assume $d\le n$. For a partition $\lambda\vdash d$, let $s_\lambda$ and $m_\lambda$ be the Schur and monomial symmetric polynomials in $z_1,\ldots,z_n$, let $\lambda'$ be the conjugate partition, and let $f^\lambda=f^{\lambda'}$ be the number of standard Young tableaux of shape $\lambda$. The immanant of $A$ attached to the irreducible character $\chi^\lambda$ of the symmetric group $S_d$ is
\[
 \Imm_\lambda(A)=\sum_{\sigma\in S_d}\chi^\lambda(\sigma)\nlprod_{i=1}^dA_{i\sigma(i)},
\]
so that $\Imm_{(1^d)}=\det$ and $\Imm_{(d)}=\per$. We write a symmetric polynomial $p$, homogeneous of degree $d$, as $p=\nlsum_{\lambda\vdash d}a_\lambda s_\lambda$. Since $s_{(d)}$ is the only Schur polynomial containing the monomial $z_1^d$, the coefficient $a_{(d)}$ is the coefficient of $z_1^d$ in $p$. Borcea and Br\"and\'en observed that the Schur coefficients of $f_A$ are immanants~\citep[before Problem~35]{aim2007polyaschurlax}:
\[
 f_A=\nlsum_{\lambda\vdash d}\Imm_{\lambda'}(A)\,s_\lambda.
\]
In particular $a_{(d)}=\det A$ and $a_{(1^d)}=\per A$ for $p=f_A$. For positive semidefinite $A$, Schur's inequality $\Imm_\lambda(A)\ge f^\lambda\det A$ and Lieb's permanental dominance conjecture $\Imm_\lambda(A)\le f^\lambda\per A$, which is known for $d\le13$~\citep[Section~3]{wanless2022lieb}, therefore become
\begin{equation}\label{eq:immanant-bounds}
 f^\lambda a_{(d)}\le a_\lambda\le f^\lambda a_{(1^d)},\qquad \lambda\vdash d,
\end{equation}
the upper bound for $d\le13$. Since $a_{(d)}=\det A\ge0$, every $f_A$ is Schur positive.

A polynomial with real coefficients is \emph{real stable} if it does not vanish when every variable has positive imaginary part; every nonzero $f_A$ is real stable, symmetric, and homogeneous~\citep{branden08}. Let $p$ be a symmetric, real stable, homogeneous polynomial of degree $d$ in $n$ variables. Problem~35 of~\citet{aim2007polyaschurlax} asks whether $p=f_A$ for some positive semidefinite matrix $A$ of order $d$. Problems~36--38 assume in addition that $d\le n$ and that $p$ has a positive coefficient, so that all its coefficients are nonnegative. They ask whether $p$ is Schur positive (Problem~36), and whether $p$ satisfies the lower bound in~\eqref{eq:immanant-bounds} (Problem~37, a stable Schur inequality) or the upper bound (Problem~38, called there the stable permanent-on-top inequality). Because $a_{(d)}\ge0$, the lower bound implies Schur positivity, so a counterexample to Problem~37 says something new only if it is Schur positive. All four counterexamples have $d=n=5$, the setting of Corollary~\ref{cor:target-intro}.

The first counterexample is a product of linear forms,
\[
 g(z)=\nlprod_{i=1}^5\Bigl(z_i+5\nlsum_{j\ne i}z_j\Bigr).
\]
Each factor has positive coefficients, so it has positive imaginary part when every $z_j$ does, and $g$ is real stable. It is symmetric and homogeneous of degree five in five variables, and all its coefficients are positive. It is also a mixed determinant of diagonal matrices. Let $D_j$ agree with $5I_5$ except that its $j$th diagonal entry is $1$. For diagonal arguments the mixed determinant factors as $\nlprod_k\nlsum_jz_j(D_j)_{kk}$, so $g(z)=\eta(z_1D_1,\ldots,z_5D_5)$. Proposition~\ref{prop:psd-tuple} therefore covers $g$; the next theorem shows that Corollary~\ref{cor:target-intro} does not.

\begin{dashedbox}[supporting]
\begin{theorem}\label{thm:aim-g}
The polynomial $g$ has the Schur expansion
\[
 g=625s_{(5)}+4500s_{(4,1)}+7125s_{(3,2)}+8150s_{(3,1,1)}+7525s_{(2,2,1)}+6580s_{(2,1,1,1)}+1281s_{(1^5)}.
\]
It violates the upper bound $a_\lambda\le f^\lambda a_{(1^5)}$ at $\lambda=(3,2),(3,1,1),(2,2,1),(2,1,1,1)$, so it answers Problem~38 negatively. There is no positive semidefinite $A\in\HH_5(\C)$ with $g=f_A$, so it also answers Problem~35 negatively.
\end{theorem}
\end{dashedbox}
\begin{proof}
In the monomial basis we have
\begin{align*}
 g={}&625m_{(5)}+5125m_{(4,1)}+12250m_{(3,2)}+24900m_{(3,1,1)}\\
 &+39550m_{(2,2,1)}+81580m_{(2,1,1,1)}+168376m_{(1^5)},
\end{align*}
and inverting the Kostka matrix, which expresses each $s_\lambda$ in the monomial basis, over $\mathbb Z$ we get the Schur expansion. The hook-length formula gives $f^\lambda=5,6,5,4$ for the four partitions, while $a_{(1^5)}=1281$, and
\[
 7125>6405,\qquad 8150>7686,\qquad 7525>6405,\qquad 6580>5124.
\]
If $g=f_A$ for a positive semidefinite matrix $A$ of order five, then the upper bound in~\eqref{eq:immanant-bounds} holds, since $5\le13$, and these inequalities contradict it.
\end{proof}

The polynomial $g$ is Schur positive and satisfies the lower bound in~\eqref{eq:immanant-bounds}, so for Problems~36 and~37 we need other witnesses, and we take both from one family of positive definite pencils. Let
\begingroup
\scriptsize
\[
T_1=
\begin{pmatrix}
38&19&19&2&36\\
19&38&2&19&36\\
19&2&38&19&36\\
2&19&19&38&36\\
36&36&36&36&72
\end{pmatrix},\quad
T_2=
\begin{pmatrix}
38&19&19&17&21\\
19&38&17&19&21\\
19&17&38&34&36\\
17&19&34&38&36\\
21&21&36&36&42
\end{pmatrix},\quad
T_3=
\begin{pmatrix}
38&34&19&17&36\\
34&38&17&19&36\\
19&17&38&19&21\\
17&19&19&38&21\\
36&36&21&21&42
\end{pmatrix},
\]
\[
T_4=
\begin{pmatrix}
38&4&19&2&21\\
4&8&2&4&6\\
19&2&38&4&21\\
2&4&4&8&6\\
21&6&21&6&42
\end{pmatrix},\qquad
T_5=
\begin{pmatrix}
8&4&4&2&6\\
4&38&2&19&21\\
4&2&8&4&6\\
2&19&4&38&21\\
6&21&6&21&42
\end{pmatrix},
\]
\endgroup
and put $h_5(z)=\frac1{648}\det\bigl(\nlsum_{i=1}^5z_iT_i\bigr)$.

\begin{dashedbox}[supporting]
\begin{theorem}\label{thm:aim-h5}
The polynomial $h_5$ is symmetric, homogeneous of degree five in five variables, and real stable, and all $126$ of its monomial coefficients are positive. It is not Schur positive, since $[s_{(1^5)}]h_5=-2972$. Hence it answers Problem~36 negatively.
\end{theorem}
\end{dashedbox}
\begin{proof}
The leading principal minors of the $T_i$ are
\[
\begin{array}{c|rrrrr}
 &\Delta_1&\Delta_2&\Delta_3&\Delta_4&\Delta_5\\ \hline
T_1&38&1083&28728&202176&1119744\\
T_2&38&1083&28728&202176&1119744\\
T_3&38&288&8208&202176&1119744\\
T_4&38&288&8208&46656&1119744\\
T_5&8&288&1728&46656&1119744,
\end{array}
\]
all positive, so each $T_i$ is positive definite by Sylvester's criterion. If every $z_i$ has positive imaginary part, the imaginary part $\nlsum_i(\operatorname{Im}z_i)T_i$ of the pencil is positive definite, so the pencil is nonsingular and $h_5(z)\ne0$. Expanding the determinant exactly we get
\begin{align*}
 h_5={}&1728m_{(5)}+19440m_{(4,1)}+66555m_{(3,2)}+140910m_{(3,1,1)}\\
 &+250440m_{(2,2,1)}+547405m_{(2,1,1,1)}+1182860m_{(1^5)},
\end{align*}
so $h_5$ is symmetric with positive coefficients, and inverting the Kostka matrix we get
\[
 h_5=1728s_{(5)}+17712s_{(4,1)}+47115s_{(3,2)}+56643s_{(3,1,1)}+62415s_{(2,2,1)}+56437s_{(2,1,1,1)}-2972s_{(1^5)}.
\]
We can check the last coefficient by hand: the last column of the inverse Kostka matrix is $(1,-2,-2,3,3,-4,1)$, so
\[
 [s_{(1^5)}]h_5=1728-2(19440)-2(66555)+3(140910)+3(250440)-4(547405)+1182860=-2972.\qedhere
\]
\end{proof}

\begin{rmk}[Origin of the pencil]\label{rmk:specht}
The matrices $T_i$ come from the Specht module $S^{(3,2)}$ of $S_5$~\citep[Appendix~C]{Sra2026CountExMachina}. Let $V\cong S^{(3,2)}$ be the kernel of the map from the permutation module on $2$-subsets of $[5]$ to $\R^5$ that sends $e_{ab}$ to $e_a+e_b$, and let $P_i$ be half the sum of the transpositions of $S_5$ that fix $i$, acting on $V$; it is an orthogonal projection of rank three. In a suitable basis of $V$ with Gram matrix $G$, one has $T_i=2G(I+5P_i)$. So $h_5$ is a multiple of $\det(\nlsum_iz_i(I+tP_i))$ at $t=5$, whose coefficient of $s_{(1^5)}$ is negative once $t$ is large.
\end{rmk}

For Problem~37 we take the member of the same family at $t=4$. Let
\[
G=\begin{pmatrix}
4&2&2&1&3\\
2&4&1&2&3\\
2&1&4&2&3\\
1&2&2&4&3\\
3&3&3&3&6
\end{pmatrix}
\]
be the Gram matrix of Remark~\ref{rmk:specht}. The matrices $\widetilde T_i=(4T_i+2G)/5=2G(I+4P_i)$ are again symmetric integer matrices. Put $h_4(z)=\frac1{5184}\det\bigl(\nlsum_{i=1}^5z_i\widetilde T_i\bigr)$.

\begin{dashedbox}[supporting]
\begin{theorem}\label{thm:aim-h4}
The polynomial $h_4$ satisfies the hypotheses of Problems~36--38 and is Schur positive, with
\[
 h_4=125s_{(5)}+1100s_{(4,1)}+2545s_{(3,2)}+3110s_{(3,1,1)}+3321s_{(2,2,1)}+2972s_{(2,1,1,1)}+69s_{(1^5)}.
\]
At $\lambda=(1^5)$, where $f^{(1^5)}=1$, it violates the lower bound in~\eqref{eq:immanant-bounds}, since $a_{(1^5)}=69<125=a_{(5)}$. Hence it answers Problem~37 negatively.
\end{theorem}
\end{dashedbox}
\begin{proof}
The leading principal minors of the $\widetilde T_i$ are
\[
\begin{array}{c|rrrrr}
 &\Delta_1&\Delta_2&\Delta_3&\Delta_4&\Delta_5\\ \hline
\widetilde T_1&32&768&17280&118800&648000\\
\widetilde T_2&32&768&17280&118800&648000\\
\widetilde T_3&32&240&5760&118800&648000\\
\widetilde T_4&32&240&5760&32400&648000\\
\widetilde T_5&8&240&1440&32400&648000,
\end{array}
\]
all positive, so each $\widetilde T_i$ is positive definite, and $h_4$ is real stable as before. Expanding exactly we get
\begin{align*}
 h_4={}&125m_{(5)}+1225m_{(4,1)}+3770m_{(3,2)}+7980m_{(3,1,1)}\\
 &+13846m_{(2,2,1)}+30004m_{(2,1,1,1)}+64472m_{(1^5)},
\end{align*}
so $h_4$ is symmetric and homogeneous of degree five in five variables with positive coefficients, and inverting the Kostka matrix we get the Schur expansion, in which the lower bound fails only at $\lambda=(1^5)$.
\end{proof}

By~\eqref{eq:immanant-bounds} every $f_A$ is Schur positive, so Theorem~\ref{thm:aim-h5} also shows that $h_5\ne f_A$ for every positive semidefinite $A$; this second answer to Problem~35 uses Schur's inequality in place of permanental dominance at order five. Since $h_5$ is the determinant of a positive definite pencil, its cone is spectrahedral, and we conclude that the family of Corollary~\ref{cor:target-intro} misses even symmetric polynomials with definite determinantal representations.

\section{A cubic whose cone needs a multiplier}\label{app:cubic}
Theorem~\ref{thm:main} represents each cone by a pencil whose determinant carries a polynomial multiplier. In this appendix we exhibit an explicit cubic $P$, a linear section of a leading coefficient $Q_{A_0}$ from Section~\ref{sec:localization}, for which the multiplier is indispensable, and we show why an obstruction modeled on~\citet{openai2026nonspectrahedral} does not apply to its cone.

\subsection{A six-variable cubic inside a rank-deficient localization}\label{subapp:cubic}
We set
\begin{equation}\label{eq:V-frame}
 V=\begin{pmatrix}
 0&0&0&1&1&1\\
 0&1&1&0&0&-1\\
 1&0&-1&0&-1&0
 \end{pmatrix},\qquad A_0=V^{\mathsf T}V.
\end{equation}
Then $A_0\succeq0$ has rank three, and we consider $Q(G)=\nlsum_{|S|=3}\det A_0[S]\det G[S^c]$ for $G\in\HH_6(\R)$, which has sixteen nonzero coefficients, each equal to one. On the linear section
\begin{equation*}
 G(x,y,z,a,b,c)=
 \begin{pmatrix}
 x&0&0&0&0&0\\
 0&y&a&b&c&0\\
 0&a&y&c&b&0\\
 0&b&c&y&a&0\\
 0&c&b&a&y&0\\
 0&0&0&0&0&z
 \end{pmatrix},
\end{equation*}
we obtain
\begin{align*}
 P(x,y,z,a,b,c)
 &=4y(x+y)(z+y)-2x(a^2+c^2)-2z(b^2+c^2)\\
 &\hspace{17mm}-4y(a^2+b^2+c^2)+8abc.
\end{align*}
The cubic $P$ is the one in the principal-minor lift of the spanning-tree polynomial of the complete graph $K_4$ studied in~\citep[Section~6]{blekherman2023linear}, with coordinates renamed, and we verify the identity $Q(G(x,y,z,a,b,c))=P$ by summing the sixteen $3\times3$ principal determinants. Also, $[t^3]\eta(tA_0,G)=Q(G)$, so $P$ is a linear section of the leading coefficient $Q_{A_0}$ of Section~\ref{sec:localization}. The section contains $G=I_6$, corresponding to $(x,y,z,a,b,c)=(1,1,1,0,0,0)$.

\subsection{The Rayleigh quartic is not a sum of squares}\label{subapp:F}
For a polynomial $p$ we write $\Delta_{xz}p=p_xp_z-pp_{xz}$. Direct calculation gives $\Delta_{xz}P=4F$ for the quartic
\begin{equation*}
 F=a^2b^2+a^2c^2+b^2c^2+c^4-8abcy+2a^2y^2+2b^2y^2.
\end{equation*}
\citet{blekherman2023linear} show that this Rayleigh quartic is not a sum of squares; we give a direct argument that also determines exactly which homogeneous quadratic squares it dominates.

\begin{prop}\label{prop:dominated}
The polynomial $F$ is nonnegative and
\begin{equation}\label{eq:dominated-space}
 \left\{q\in\R[a,b,c,y]_2:
 q^2\le C F\text{ on }\R^4\text{ for some }C>0\right\}
 =\Span\{ac-by,\ ay-bc\}.
\end{equation}
Consequently, $F$ is not a sum of polynomial squares and there is no finite matrix $M(a,b,c,y)$ of homogeneous quadratic polynomials such that
\begin{equation}\label{eq:impossible-norm}
 F=\norm{M}_{\op}^2.
\end{equation}
\end{prop}
\begin{proof}
We set
\[
 C_0=a^2b^2+a^2y^2+b^2y^2+c^4-4abcy.
\]
By the AM--GM inequality $C_0\ge0$, and expanding we find $F=C_0+(ac-by)^2+(ay-bc)^2$; thus $F\ge0$, and each generator on the right of~\eqref{eq:dominated-space} has square at most a multiple of $F$.

Conversely, suppose $q^2\le CF$. The $a$-, $b$-, and $y$-axes are zeros of $F$, so the coefficients of $a^2,b^2,y^2$ in $q$ vanish, and comparing the four sign characters at the further zeros $(a,b,c,y)=(\sigma,\tau,\sigma\tau,1)$ with $\sigma,\tau\in\{\pm1\}$, we get $q=\alpha(ab-cy)+\beta(ac-by)+\gamma(ay-bc)$. On the line $(a,b,c,y)=(0,0,s,1)$ we have $F=s^4$ and $q=-\alpha s$, and $s\to0$ in the domination inequality forces $\alpha=0$, proving~\eqref{eq:dominated-space}.

At $(a,b,c,y)=(0,0,1,0)$, both remaining generators vanish but $F=1$. In a sum-of-squares expression for a homogeneous quartic we may take each nonzero summand homogeneous quadratic, since scaling and the absence of cancellation between squares eliminate other degrees; each such summand would be dominated by $F$ and hence would vanish at this point, a contradiction. Likewise, every entry of a matrix in~\eqref{eq:impossible-norm} would have its square bounded by $F$, so every entry would vanish there, again a contradiction.
\end{proof}

We let $e_y=(0,1,0,0,0,0)$ in the coordinate order $(x,y,z,a,b,c)$; it determines the same component as $(1,1,1,0,0,0)$, since $P(\tau,1,\tau,0,0,0)=4(1+\tau)^2>0$ for $0\le\tau\le1$. The coordinate directions $e_x,e_z$ belong to the closed cone; for example $P(e_x+te_y)=4t^2(t+1)$, and similarly for $e_z$.

Since $\Delta_{xz}P=4F$ is not a sum of squares, the Rayleigh criterion of \citet[Corollary~4.3]{kummer2015hyperbolic} shows that no positive power of $P$ has a definite symmetric determinantal representation for this component; we may allow boundary directions by taking limits from the open cone and using closedness of the fixed-degree sum-of-squares cone. The criterion obstructs \emph{pure powers} only; it does not exclude a determinant equal to $P$ times other factors whose cones contain the cone of $P$.

An obstruction to the cone itself, in the style of~\citet{openai2026nonspectrahedral}, would need the implication
\[
 \text{the cone of }P\text{ is spectrahedral}
 \quad\Longrightarrow\quad
 F=\norm{M}_{\op}^2\text{ for a quadratic polynomial matrix }M,
\]
whose conclusion Proposition~\ref{prop:dominated} forbids. With the pencil of Appendix~\ref{subapp:17} we show that this implication fails, so Proposition~\ref{prop:dominated} does not obstruct spectrahedrality.

\subsection{A \texorpdfstring{$17\times17$}{17x17} pencil for the cubic}\label{subapp:17}
We make the multiplier explicit with a pencil whose identities we checked in exact arithmetic over $\Q(\sqrt2,\sqrt3)$. We use the invertible linear substitution
\begin{align*}
 t&=y+\frac{x+z}{3},&u&=-\frac{x+z}{6},&v&=\frac{x-z}{2\sqrt3},\\
 r&=\frac a{\sqrt3},&s&=\frac b{\sqrt3},&w&=\frac c{\sqrt3},
\end{align*}
and put $N=u^2+v^2+r^2+s^2+w^2$. For $\xi=(u,v,\alpha,\beta,\gamma)^{\mathsf T}$ we define
\begin{align*}
 q_C(\xi)&=u^3-3uv^2+\frac32u(\alpha^2+\beta^2-2\gamma^2)\\
 &\hspace{10mm}+\frac{3\sqrt3}{2}v(\alpha^2-\beta^2)
       +3\sqrt3\alpha\beta\gamma.
\end{align*}
We also use the matrix
\begin{equation*}
 U=\frac{\sqrt2}{4}
 \begin{pmatrix}
 2&0&0&0&2\\
 0&0&1&\sqrt3&0\\
 2&0&0&0&-2\\
 0&0&-1&\sqrt3&0\\
 0&2&0&0&0\\
 0&0&\sqrt3&1&0\\
 0&-2&0&0&0\\
 0&0&\sqrt3&-1&0
 \end{pmatrix}.
\end{equation*}
We define the $4\times4$ matrices
\begin{align*}
 A(\xi)&=\begin{pmatrix}
 -\gamma&\alpha&u&0\\
 \alpha&\gamma&0&u\\
 u&0&\gamma&-\alpha\\
 0&u&-\alpha&-\gamma
 \end{pmatrix},\\[1mm]
 B(\xi)&=\begin{pmatrix}
 0&-\beta&-v&0\\
 \beta&0&0&-v\\
 v&0&0&\beta\\
 0&v&-\beta&0
 \end{pmatrix},\qquad
 \Gamma(\xi)=\begin{pmatrix}A(\xi)&-B(\xi)\\B(\xi)&A(\xi)\end{pmatrix}.
\end{align*}
Direct multiplication gives
\begin{align}
 U^{\mathsf T}U&=I_5,&\Gamma(\xi)^{\mathsf T}&=\Gamma(\xi),\notag\\
 \Gamma(\xi)^2&=\norm{\xi}^2I_8,&\tr\Gamma(\xi)&=0,\label{eq:clifford-identities}\\
 (U\xi)^{\mathsf T}\Gamma(\xi)(U\xi)&=q_C(\xi).&&\notag
\end{align}
In particular, $|q_C(\xi)|\le\norm{\xi}^3$.

We take the two vectors $\xi_1=(u,v,s,w,r)^{\mathsf T}$ and $\xi_2=(u,v,w,r,s)^{\mathsf T}$, both of squared norm $N$, and direct expansion gives
\begin{equation*}
 h(t,u,v,r,s,w):=t^3-3Nt+q_C(\xi_1)+q_C(\xi_2)=\frac P4.
\end{equation*}
We now define the $17\times17$ real symmetric linear pencil
\begin{equation}\label{eq:L17}
 L_P=\begin{pmatrix}
 t&(U\xi_1)^{\mathsf T}&(U\xi_2)^{\mathsf T}\\
 U\xi_1&tI_8+\Gamma(\xi_1)&0\\
 U\xi_2&0&tI_8+\Gamma(\xi_2)
 \end{pmatrix}.
\end{equation}

\begin{prop}\label{prop:L17}
The pencil~\eqref{eq:L17} satisfies
\begin{equation}\label{eq:L17-result}
 L_P(e_y)=I_{17},\qquad
 \cone(P,e_y)=\{L_P\succeq0\},\qquad
 \det L_P=\frac P4\left(\frac{P_y}{12}\right)^7.
\end{equation}
\end{prop}
\begin{proof}
We write $q=\tfrac12(q_C(\xi_1)+q_C(\xi_2))$. Then $|q|\le N^{3/2}$, and $h=t^3-3Nt+2q$ has only real roots: its values at the two critical points satisfy $h(-\sqrt N)\ge0$ and $h(\sqrt N)\le0$, and when $N=0$ it is $t^3$. Its largest root $\rho$ satisfies $\rho\ge\sqrt N$, and $h$ is strictly increasing for $t>\sqrt N$.

The lower blocks of~\eqref{eq:L17} are positive definite exactly when $t>\sqrt N$, and on that region $(tI_8+\Gamma(\xi_i))^{-1}=(tI_8-\Gamma(\xi_i))/(t^2-N)$. Using~\eqref{eq:clifford-identities}, we compute the scalar Schur complement
\[
 t-\frac{2tN-q_C(\xi_1)-q_C(\xi_2)}{t^2-N}
       =\frac{h}{t^2-N}.
\]
Thus $L_P\succ0$ exactly when $t>\rho$, which is the open hyperbolicity cone of $h$ in the $t$-direction. Increasing $y$ with $x,z,a,b,c$ fixed increases $t$ alone, so this direction is $e_y$ in the original coordinates, and since $L_P(e_y)=I$, closure gives the cone equality.

The square and trace identities imply
$\det(tI_8+\Gamma(\xi_i))=(t^2-N)^4$, so the Schur determinant formula gives $\det L_P=h(t^2-N)^7$, initially on a dense open set and therefore polynomially. Finally, expanding we find $t^2-N=P_y/12$, which proves~\eqref{eq:L17-result}.
\end{proof}

The multiplier $(P_y/12)^7$ in~\eqref{eq:L17-result} escapes the pure-power criterion of Appendix~\ref{subapp:F}.

\subsection{The full cubic and its zero border}\label{app:full-cubic}
The section $P$ comes from the full $21$-variable cubic $Q$ of Appendix~\ref{subapp:cubic}, which also has a compact form. We reverse the columns of $V$ in~\eqref{eq:V-frame} to obtain
\[
 W=\begin{pmatrix}
 1&1&1&0&0&0\\
 -1&0&0&1&1&0\\
 0&-1&0&-1&0&1
 \end{pmatrix},
\]
with columns $w_1,\ldots,w_6$, write $\mu_i=G_{ii}$ and $g_{ij}=G_{ij}$ for $i<j$, and define
\begin{align*}
 M&=W\diag(\mu_1,\ldots,\mu_6)W^{\mathsf T},\\
 B(g)&=\nlsum_{i<j}g_{ij}^2(w_i\times w_j)(w_i\times w_j)^{\mathsf T},\\
 \tau(g)&=\nlsum_{\{i,j,k\}\in\mathcal T}g_{ij}g_{ik}g_{jk},\\
 \mathcal T&=\binom{[6]}3\setminus\{124,135,236,456\}.
\end{align*}
The map $\mu\mapsto M$ is an isomorphism from $\R^6$ to $\HH_3(\R)$: the three off-diagonal entries of $M$ recover $\mu_1,\mu_2,\mu_4$, and the diagonal entries recover $\mu_3,\mu_5,\mu_6$. Expanding the principal minors we get $Q(M,g)=\det M-\tr(MB(g))+2\tau(g)$. Its zero-border extension $Q^\sharp$ of~\eqref{eq:zero-border}, obtained with $\widehat A_0=\diag(A_0,0)$, has $28$ variables; the factor $\tfrac12$ on off-diagonal coordinates in the gradient convention~\eqref{eq:matrix-gradient} matters there. An obstruction to the spectrahedrality of $Q^\sharp$ along the lines of Appendix~\ref{subapp:F} would have to show that no fixed linear maps $\Phi_0$ and $R$ satisfy the norm identity~\eqref{eq:resolvent-identity} throughout the interior of the cone of $Q$. In Proposition~\ref{prop:compression} and Corollary~\ref{cor:zero-border} we constructed exactly such maps, so the cones of the full cubic and of its zero-border extension are both spectrahedral.

\end{document}